\documentclass[11pt]{amsart}

\usepackage[square,compress,comma, numbers,sort]{natbib}
\usepackage[colorlinks=true, citecolor=red, linkcolor=blue]{hyperref}
\usepackage{amsfonts,mathtools}

\allowdisplaybreaks[4]

\usepackage{amsmath}
\usepackage{bbm}
\usepackage{amssymb}
\usepackage{mathtools}

\usepackage{color}

\def\Ee#1{{\textcolor{black}{#1}}}

\def\II{\mathcal{I}}
\def\JJ{\mathcal{J}}
\def\KK{\mathcal{K}}

\definecolor{c20}{rgb}{0.,0.7,0.}
\definecolor{c30}{rgb}{0.,0.,1.}
\definecolor{c40}{rgb}{1,0.1,0.7}
\definecolor{c50}{rgb}{1,0,0}
\definecolor{c60}{rgb}{1,0.9,0.1}
\definecolor{krzys}{rgb}{0.8,0.1,1}
\definecolor{brown}{rgb}{0.5,0.38,0.3}
\def\EZ#1{{\textcolor{black}{#1}}}

\def\CE#1{{\textcolor{cyan}{#1}}}
\def\CE#1{{\textcolor{black}{#1}}}

\def\CN#1{{\textcolor{black}{#1}}}
\def\CNN#1{{\textcolor{black}{#1}}}

\def\kk#1{{\textcolor{black}{#1}}}
\def\kdd#1{{\textcolor{black}{#1}}}

\def\cE#1{{\textcolor{black}{#1}}}

\def\k1#1{{\textcolor{black}{#1}}}

\def\x{\vk{x}}

\newcommand{\abs}[1]{\left\lvert #1 \right\rvert}

\newcommand{\sprod}[1]{\langle#1\rangle}

\newcommand{\E}[1]{\mathbb{E}\left\{ #1\right\}}

\newcommand{\pk}[1]{\mathbb{P} \left\{ #1 \right \} }

\newcommand{\R}{\mathbb{R}}

\newcommand{\inr}{\in \R}

\newcommand{\ldot}{,\ldots,}

\newcommand{\limit}[1]{\lim_{#1 \to \infty}}

\newcommand{\BQN}{\begin{eqnarray}}
\newcommand{\EQN}{\end{eqnarray}}
\newcommand{\BQNY}{\begin{eqnarray*}}
\newcommand{\EQNY}{\end{eqnarray*}}

\newcommand{\BS}{\begin{sat}}
\newcommand{\ES}{\end{sat}}
\newcommand{\BT}{\begin{theo}}
\newcommand{\ET}{\end{theo}}
\newcommand{\BK}{\begin{korr}}
\newcommand{\EK}{\end{korr}}

\newcommand{\BD}{\begin{de}}
\newcommand{\ED}{\end{de}}

\newcommand{\BIT}{\begin{itemize}}
\newcommand{\EIT}{\end{itemize}}
\newcommand{\BDI}{\begin{description}}
\newcommand{\EDI}{\end{description}}

\newcommand{\BRM}{\begin{remark}}
\newcommand{\ERM}{\end{remark}}

\newcommand{\BEL}{\begin{lem}}
\newcommand{\EEL}{\end{lem}}

\newtheorem{theo}{Theorem}[section]
\newtheorem{sat}[theo]{Proposition}
\newtheorem{de}[theo]{Definition}
\newtheorem{lem}[theo]{Lemma}

\newtheorem{korr}[theo]{Corollary}
\newtheorem{remark}[theo]{Remark}

\newcommand{\nelem}[1]{{Lemma \ref{#1}}}

\newcommand{\netheo}[1]{{Theorem \ref{#1}}}

\newcommand{\prooftheo}[1]{ \textsc{\bf Proof of Theorem} \ref{#1}:}

\newcommand{\prooflem}[1]{\textsc{\bf Proof of Lemma} \ref{#1}:}

\newcommand{\COM}[1]{}

\def\td{\text{\rm d}}
\def\td{d}
\def\IF{\infty}

\newcommand{\QED}{\hfill $\Box$}

\newcommand{\kb}[1]{\boldsymbol{#1}}
\newcommand{\vk}[1]{\kb{#1}}

\def\IF{\infty}

\def\bqny#1{\begin{eqnarray*} #1 \end{eqnarray*}}
\def\bqn#1{\begin{eqnarray} #1 \end{eqnarray}}

\begin{document}


\title{Pandemic-type Failures in   Multivariate Brownian Risk Models}

\author{Krzysztof D\c{e}bicki}
\address{Krzysztof D\c{e}bicki, Mathematical Institute, University of Wroc\l aw, pl. Grunwaldzki 2/4, 50-384 Wroc\l aw, Poland}
\email{Krzysztof.Debicki@math.uni.wroc.pl}

\author{Enkelejd Hashorva}
\address{Enkelejd Hashorva, Department of Actuarial Science, 
University of Lausanne,\\
UNIL-Dorigny, 1015 Lausanne, Switzerland
}
\email{Enkelejd.Hashorva@unil.ch}

\author{Nikolai Kriukov}
\address{Nikolai Kriukov, Department of Actuarial Science, 
University of Lausanne,\\
UNIL-Dorigny, 1015 Lausanne, Switzerland
}
\email{Nikolai.Kriukov@unil.ch}

\bigskip

\date{\today}
\maketitle

{\bf Abstract:} Modelling of multiple simultaneous failures in insurance, finance and other areas of applied probability
is important especially from the point of view of pandemic-type events.
A benchmark limiting model for the analysis of multiple failures is the  classical $d$-dimensional Brownian risk model (Brm), see \cite{delsing2018asymptotics}.
From both theoretical and practical point of view, of interest is the calculation of
the probability of multiple simultaneous failures in a given time horizon.
The  main findings of this contribution concern the approximation of  the probability that at least $k$ out of $d$ components of
Brm fail simultaneously.  We derive both sharp bounds and asymptotic approximations of the probability of interest
for the finite and the infinite time horizon.  Our results extend previous findings of \cite{mi:18,Rolski17}.

{\bf Key Words:} Multivariate Brownian risk model; probability of multiple simultaneous failures; simultaneous ruin probability;
 failure  time; exact asymptotics; pandemic-type events.

{\bf AMS Classification:} Primary 60G15; secondary 60G70\\

\section{Introduction}
In this paper we are interested in the probabilistic aspects of multiple simultaneous failures typically occurring due to pandemic-type events.
A key benchmark risk model  considered here is the $d$-dimensional Brownian risk model (Brm)
$$ \vk R(t, \vk u)\kdd{=(R_1(t,u_1),\ldots,R_d(t,u_d))^\top}= \vk u   \k1{+ \vk c t- \vk W(t)}, \quad t\ge 0,$$
where
$\vk c= (c_1 \ldot c_d)^{\kdd{\top}}, \vk u=(u_1 \ldot u_d)^\top $ are vectors in $\R^d$
 and
$$\vk W(t)= \Gamma  \vk B(t), \quad t\in \R ,$$
with $\Gamma$ a $d\times d$ real-valued \CN{non-singular} matrix and $\vk B(t)=(B_1(t) \ldot B_d(t))^\top , t \in \R$ a $d$-dimensional
Brownian motion with independent components which are standard Brownian motions.

By  bold symbols we denote column vectors, operations with vectors are meant component-wise and $ a \vk x= (a x_1 \ldot a x_d)^\top$ for any scalar $a\inr$ and any $\vk x\inr^d$.

 Indeed, Brm is a natural limiting model in many statistical applications.  Moreover, as shown in  \cite{delsing2018asymptotics} such a risk
  model  appears naturally in insurance applications. Being a limiting model,
  Brm can be used as a benchmark for various untractable  models.
  Given the fundamental role of Brownian motion in applied probability and statistics, it is also of theoretical interest to study failure events arising from this model. Specifically,
  in this contribution we are interested in the behaviour of the probability of multiple simultaneous
  failures occurring in a  given time horizon  $[S,T] \subset [0, \IF]$.

In our settings failures can be defined in various ways.
Let us consider first the failure of a given component of our risk model. Namely, we say that the $i$th component  of our Brm has a failure (or  ruin occurs)
if $ R_i(t, \kdd{u_i})\CNN{=u_i\k1{+c_i t-W_i(t)}}< 0$ for some $t\in [S,T]$.
The extreme case of a catastrophic  event   is when  $d$ multiple simultaneous failures  occurs.
Typically, for pandemic-type events there are at least $k$ components of the model with simultaneous failures and $k$ is large with the extreme case $k=d$.
In mathematical notation, for given positive integer $k \le d$ of interest is the calculation of the following probability
\bqny{
	\psi_k(S,T,\vk u)&=& \pk{\exists t\in[S,T],\,\exists \II\subset\{1\ldot d\},|\II|=k: \cap_{i\in \II} \{ R_i(t, \kdd{ u_i}) < 0\} } \\
	&=& \pk{\exists t\in[S,T],\,\exists \II\subset\{1\ldot d\},|\II|=k: \kk{\cap_{i\in \II}} \, \CNN{\{}W_i(t)- c_i t>u_i\CNN{\}}},
}
where $|\II|$ denotes the cardinality of the set $\II$. If $T$ is finite,
by the self-similarity
property of the Brownian motion $\psi_k(S,T,\vk u)$  can be   derived from the case $T=1$,
whereas $T=\IF$ has to be treated separately.

There are no results in the literature investigating $\psi_k(S,T,\vk u)  $ for general $k$. The particular case $k=d$, for which $\psi_d(S,T,\vk u)  $ coincides with
the simultaneous ruin probability  has been studies in different contexts, see e.g.,  \cite{Mandjes07, PistAvram1, PistAvram2,Rolski17, MR3776537, foss2017two, borovkov17, PalmBorov,ji2018cumulative, MR3102482,MR3493174, cl2017}. The case $d=2$ of Brm has been recently investigated in  \cite{mi:18}.

Although the probability of multiple simultaneous failures seems very difficult to compute,
our first result below, motivated by \cite{KWW}[Thm 1.1], shows that \kk{$\psi_k(S,T,\vk u)$} can be bounded by  the  multivariate Gaussian survival probability, namely by
$$ p_T(\vk u)= \pk{\Ee{(W_1(T)- c_1 T \ldot W_d(T)- c_d T)} \in \vk E(\vk u)},$$
\kk{
where
\bqn{\vk E(\vk u)\CNN{=\bigcup_{\substack{\II\subset\{1\ldot d\}\\|\II|=k}}\vk{E_\II(\vk u)}}=\bigcup_{\substack{\II\subset\{1\ldot d\}\\|\II|=k}}\Ee{\{\vk x\in\R^d:\forall i\in \II:~\vk x_i\geq\vk u_i\}}\label{def.E_I}.
}
}
\Ee{When $u\to \IF$ we can approximate   $p_T(\vk u)$ utilising Laplace asymptotic method, see e.g., \cite{MR3399144},
whereas for small and moderate values of $u$ it can be calculated or simulated with sufficient accuracy.
Our next \k1{result} gives  bounds for
$ \psi_k(\CE{S},T, \vk u)$ in terms of  $p_T(\vk u)$. }

\begin{theo}\label{D}
	If  the matrix $\Gamma$ is non-singular,  then for any positive integer $k \le d$, all constants $ 0\leq S < T< \IF$ and all $\vk c,\vk u\in\R^d$
	\bqn{
		p_T(\vk u) \le \psi_k(\CE{S},T, \vk u)
		\le
		K p_T(\vk u),
		\label{eqB1}
	}
	where $K= 1/\min_{\substack{\II\subset\{1\ldot d\},|\II|=k}}\pk{\forall_{i\in \II}: W_i(T)> \max(0,c_i T) }>0.$
\end{theo}

The bounds in \eqref{eqB1} \Ee{indicate} that it might be possible to  derive an 	approximations of $\psi_k(\CE{S},T, \vk u)$ for large threshold $\vk u$, which has been already shown for $k=d=2$ in \cite{mi:18}.  In this paper we consider the general case $k \le  d, d> 2$ discussing both the finite time interval (i.e., $T=1$)
and the infinite time horizon case with $T=\IF$ extending the results of \cite{Rolski17} where $d=k$ is considered.

In Section 2 we explain the main ideas that lead to the  approximation of $\psi_k(\CE{S},T, \vk u)$.
   Section 3 discusses some interesting special cases, whereas the proofs are postponed to Section 4. Some technical calculations are displayed in  Section \ref{APP}.

\section{Main Results}

In this section $\vk W(t),t\ge 0$ is as in the Introduction and for a given positive integer $k\le d$ we shall investigate the  approximation of $ \psi_k(S,T, \vk u)$
where we fix $\vk u = \vk a u$,
    with $\vk a$ in $\R^d\setminus(-\infty,0]^d$ and $u$ is sufficiently large.

Let hereafter $\II$ denote a non-empty index set of $\{1 \ldot d\}$. For a given vector, say $\vk x\inr^d$ we shall write $\vk x_\II$ to denote a subvector of $\vk x$ obtained by dropping its components not in $\II$. Set next
 $$ \psi_{\II}(S,T, \vk a_\II u)= \pk{ \exists {t\in [S,T]}: A_\II(t)},$$
with
\bqn{\label{AI}
	 A_\II (t)=\{ \vk W(t)- \vk c t\in\vk E_\II(\vk a u) \}=\{ \forall i\in\II:\ ~W_i(t)- c_it\geq a_i u \},
}
\kdd{where $E_\II(\vk a u)$ was defined in (\ref{def.E_I}).}	
In vector notation for any $u\inr$ 
\bqny{
\psi_k(S,T,\vk a u) &=&
\pk{\exists t\in[S,T]:~\bigcup\limits_{\substack{\II\subset\{1\ldot d\}\\|\II|=k}}A_{\II}(t)}= \pk{\bigcup\limits_{\substack{\II\subset\{1\ldot d\}\\|\II|=k}}\left\{\exists t\in[S,T]:~A_\II(t)\right\}}.
}
\kdd{The following lower bound (by Bonferroni inequality)}
\BQN
\psi_k(S, T,\vk a u) &\ge & \sum_{\substack{\II\subset\{1\ldot d\}\\|\II|=k}} \psi_{\II} (S,T,\vk a_\II u)
-\sum_{\substack{\II,\JJ\subset\{1\ldot d\}\\|\II|=|\JJ|=k\\\II\not=\JJ}}\pk{\exists t,s\in[S,T]:~A_\II(t)\cap A_\JJ(s)}
\label{inclusoin1}
\EQN
together with the upper bound
\BQN
\psi_k(S, T,\vk a u) &\le & \sum_{\substack{\II\subset\{1\ldot d\}\\|\II|=k}} \psi_{\II} (S,T,\vk a_\II u)
\label{inclusoin2}
\EQN
are crucial
\kdd{for the derivation of  the exact asymptotics of $\psi_k(S, T,\vk a u)$ as $u\to\infty$}.
As we shall show  below, the upper bound (\ref{inclusoin2}) turns out to be exact asymptotically as $u\to \IF$.
The following theorem constitutes the main finding of this contribution.

\BT \label{theo1} Suppose that the  square $d\times d$ real-valued matrix $\Gamma$ is non-singular.
If $\vk a$ has no more than $k-1$ non-positive components,
\k1{where $k \le d$ is a positive integer},
then
for all $ 0 \leq
 S< T < \IF, \vk c\inr^d $
 \bqn{
\label{k1}
\psi_k(S, T, \vk a u) &\sim &
\sum_{\substack{\II\subset\{1\ldot d\}\\|\II|=k}} \psi_{\II} (0,T,\vk a_\II u)
, \quad u\to \IF.
}
Moreover, \eqref{k1} holds also if $T=\IF$, provided that
\CN{$\vk c$ and $\vk a+\vk c t$}
have no more than $k-1$ non-positive components \CN{for all $t\geq 0$}.
\ET

\Ee{Essentially, the above result is the claim that
	the second term in the Bonferroni
lower bound (\ref{inclusoin1}) is asymptotically negligible.} In order to prove that,  the  asymptotics of
$ \psi_{\abs{\II}} (S,T,\vk a_\II u)$ has to be derived. For the special case that $\II$ has only two elements and $S=0$,
its approximation has been obtained  in \cite{mi:18}.
\k1{The assumption in Theorem \ref{theo1} that $\vk a$ has no more than $k-1$ non-positive components excludes  the case that there exists a set $\II \subset\{1\ldot d\}, \ |\II|=k$
such that $\psi_{\II} (0,T,\vk a_\II u)$ does not tend to $0$ as $u\to\infty$,
which due to its non-rare event nature is out of interest in this contribution.}
\\
The  next result     extends the findings of \cite{mi:18} to the case $d>2$. For \kk{notational simplicity}  we consider the case $\II$ has $d$ elements and thus avoid indexing by $\II$.
Recall that in our model $\vk W(t)= \Gamma \vk B(t)$ where $\vk B(t)$ has independent standard Brownian motion components
and $\Gamma$ is a $d\times d$ non-singular real-valued matrix. Consequently $\Sigma= \Gamma \Gamma^\top $ is a positive definite matrix.

\def\b{\vk a}
\def\tilb{\tilde{\vk a}}
\cE{Hereafter $\vk 0 \in \R^d$ is the column vector with all elements equal 0.} Denote by $\Pi_\Sigma( \vk a)$ the quadratic programming problem:
$$ \text{minimise   } \vk x^\top\Sigma^{-1} \x,   \text{  for all   } \vk x\ge \vk a.$$
Its unique solution $\tilde{ \vk a}$ is such that
\begin{equation}\label{eq:IJ2}
\tilde{\vk a}_I= \vk a_I, \  \ (\Sigma_{II})^{-1} \b_{I}>\vk{0}_I, \ \ \tilde{\vk a}_J= \Sigma_{JI} (\Sigma_{II} )^{-1} \vk a_I \ge \vk a_J,	
\end{equation}
where  $\tilde{\vk a}_J$ is defined  if $J= \{1 \ldot d\} \setminus I$ is non-empty.
The index set $I$ is unique with $m=\abs{I} \ge 1$ elements, see \kk{lemma} below (or \cite{Rolski17}[Lem 2.1]) for more details.

\BEL \label{prop1} Let $\Sigma$ be a $d\times d$ positive definite matrix and let $\b \inr^d \setminus (-\IF, 0]^d $.
$\Pi_\Sigma(\b)$  has a unique solution $\tilb$ given in (\ref{eq:IJ2}) with $I$  a unique non-empty
index set  with $m\le d$ elements such that
\bqn{
\label{eq:IJ3}
\min_{\x \ge \b}\x^\top \Sigma^{-1} \x= \tilb ^\top \Sigma^{-1} \tilb
&=& \b_{I}^\top (\Sigma_{II})^{-1}\b_{I}>0,
}
\BQN \label{eq:new}
\x^\top \Sigma^{-1} \tilb= \x_F^\top (\Sigma_{FF})^{-1} \kk{\tilb_F}, \quad \forall \x\inr^d
\EQN
for any index set $F \subset\{1\ldot d \}$ containing $I$. Further if $\b= \cE{(a \ldot a)^\top , a}\in (0,\IF)$,
then $ 2 \le \abs{I} \le d$.
\EEL

In the following we set
$$\vk \lambda = \Sigma^{-1} \tilde{\vk a}.$$
In view of the above \kk{lemma}
\bqn{\label{MC}
\vk \lambda_I= ( \Sigma_{II})^{-1} \vk a_I> \vk 0_I, \  \  \vk \lambda_J \ge \vk 0_J,
}
with the convention that when $J$ is empty the indexing should be disregarded so that the last inequality above is irrelevant.

   The next  theorem extends the main result in \cite{mi:18} and further complements findings presented in Theorem \ref{theo1}
showing that the simultaneous ruin probability (i.e., $k=d$) behaves up to some constant, asymptotically as $u\to \IF$ the same as $p_T(\vk u)$. For notational simplicity and without loss of generality we consider next $T=1$.

\BT \label{theo2} If $\vk a \inr^d$ has at least one positive component and $\Gamma$ is non-singular, then for all $S\in [0,1)$
\bqn{ \psi_d (S,1, \vk a u) & \sim & C(\vk a) p_1(\vk a u)
, \quad u\to \IF,
}
where $C(\vk a)= \prod_{i\in I} \lambda_i \int_{\R^m} \pk{\exists_{t\ge 0}:\vk W_{I}(t)-t\vk a_{I}>\vk x_I}e^{ \vk \lambda^\top _I \vk x_I} \td\vk x_I \in (0,\IF)$.
\ET

\BRM\label{rem1} i)
By \nelem{inimp:Lem} below taking $T=1$ therein (hereafter $\varphi$ denotes the probability density function (pdf) of $ \kk{\Gamma}\vk B(1)$)
\bqn{ \label{p1}
p_1(\vk a u) =\pk{\vk W(1)-\vk c >u\vk a}\sim
\prod_{i\in I} \lambda_i^{-1}\pk{ \vk W_U(1)> \vk c_U \lvert
\vk W_I(1)> \vk c_I } u^{-\abs{I}}\varphi(u\tilde{\vk a}+\vk c)
}
as $u\to \IF$, where $\vk \lambda= \Sigma^{-1} \tilde{\vk a}$ and if $J= \{1 \ldot d \} \setminus I$ is non-empty,
then $U=\{j\in J:\tilde{ a}_j= a_j\}$. When $J$ is empty the \kk{conditional probability} related to $U$ above is set to 1.\\
ii) Combining \netheo{theo1} and \ref{theo2} for all $S\in [0,1)$ and all $\vk a \in \R^d$ with at least one positive component we have
 \bqn{ \psi_k(S, 1,\vk a u) \sim \sum_{\substack{\II\subset\{1\ldot d\}\\|\II|=k}} C(\vk a_\II) \psi_{\abs{\II}} (0,1,\vk a_\II u)
\sim C \pk{\forall_{i \in \II^*}: W_i(1) > ua_i+ c_i},
\quad u\to \IF
\label{lere}
}
for some $C>0$ and some $\II^*\subset \{1 \ldot d\}$ with $k$ elements.\\
iii) Comparing the results of
 \netheo{theo2} and \cite{Rolski17} we obtain $$ \limsup_{u\to \IF} \frac{ (-\ln \psi_k(S_1, 1,\vk a u))^{1/2}}{ - \ln \psi_k(S_2, \IF,\vk a u) }< \IF $$
for all  $ S_1\in [0,T], S_2\in [0,\IF)$. \\
iv)
Define  the failure time (consider for simplicity $k=d$)
for our multidimensional model by
\bqny{
\tau(u)=\inf\{t\ge 0:\vk W(t)-t\vk c>\vk a u\},\qquad u>0.
}
If $\vk a$ has at least one positive component, then for all $T>S\geq 0,x>0$
\bqn{\label{time2}
\limit{u} \pk{u^2(T-\tau(u))\geq x|\tau(u) \kk{\in [S,T]}} = {e^{-x\frac{\tilde{\vk a}^\top\Sigma^{-1}\tilde{\vk a}}{2T^2}}},
}

see the proof  in Section 4.
\ERM

\section{Examples}
In order to illustrate our findings we shall consider three examples assuming that $\Gamma \Gamma^\top$ is a positive definite correlation matrix. The first example  is dedicated to the simplest case $k=1$. In the second one we discuss $k=2$ restricting $\vk a$ to have all components equal to 1 followed then by the last example where  only the assumption  $\Gamma  \Gamma^\top $ is an equi-correlated correlation matrix is imposed.
In order to avoid additional notation we shall consider for simplicity on  $T=1$.
\cE{Below $S\in [0,1)$ is fixed.}
\\
{\bf Example 1} ($k=1$):
Suppose that $\vk a$ has all components positive. In view of \netheo{theo1} we have that
$$\psi_k(S,1, \vk a u) \sim \sum_{i=1}^d \psi_{\{i\}}(0,1, a_iu)$$
as $u\to \IF$. Note that for any positive integer $i\le d$
$$\psi_{\{i\}}(0, 1, a_i u) = \pk{ \exists_{t\in [0,1]}: B(t)- c_i t > a_i u},
$$
where $B$ is a standard Brownian motion. It follows easily that
$$\psi_k(S,1, \vk a u) \sim 2\sum_{i=1}^{d}\pk{ B(1)> a_i u+c_i}, \quad u\to \IF.$$

{\bf Example 2} ($k=2$ and $\vk a=\vk 1$): Suppose  next   $k=2$ and  $\vk a$ has all components equal 1.
 In view of Theorems \ref{theo1} and \ref{theo2} we have that
$$\psi_k(S,1,\vk 1 u)\sim\sum_{\{i,j\}\subset\{1\ldot d\}}C_{i,j}(\vk 1)\pk{\min_{k\in\{i,j\}}(W_k(1)-c_k)>u)}$$
as $u\to \IF$, where $\vk 1 \inr^d$ has all components equal to 1. By   Remark \ref{rem1} we obtain further
\bqny{
  \pk{\min_{k\in\{i,j\}}(W_k(1)-c_k)>u)}
 & \sim &\frac{u^{-2}}{(1-\rho_{i,j})^2\sqrt{2\pi(1-\rho_{i,j}^2)}}e^{-\frac{u^2}{1+\rho_{i,j}}
-\CNN{\frac{(c_i+c_j)u}{1+\rho_{i,j}}} -\frac{c_i^2-2\rho_{i,j}c_ic_j+c_j^2}{2(1-\rho^2_{i,j})}}, \quad u\to \IF.
}
Here we set $\rho_{i,j}=corr(W_i(1),W_j(1))$. Consequently, if $\rho_{i,j}>\rho_{i^*,j^*}$, then as $u\to \IF$
$$ \pk{ \min_{k\in\{i^*,j^*\}}(W_k(1)-c_k>u)}=o\left(\pk{\min_{k\in\{i,j\}}(W_k(1)-c_k>u)}\right).$$
The same holds also if $\rho_{i,j}=\rho_{i^*,j^*}$ and $c_i+c_j>c_{i^*}+c_{j^*}$. If we denote by $\tau $ the maximum of all $\rho_{i,j}$'s and by $c_*$ the maximum of $c_i+c_j$ for all $i,j$'s such that $\rho_{i,j}=\tau$, then we conclude that
$$\psi_k(S,1,\vk au)\sim\newline \sum_{i,j \in \{1\ldot d\}, \rho_{i,j}=\tau,~c_i+c_j=c_*}C_{i,j}(\vk 1)\pk{\min_{k\in\{i,j\}}(W_k(1)-c_k>u)}.
$$
 Note that in this case $C_{i,j}(\vk 1)$ does not depend on $i$ and $j$ and is equals to $$(1-\tau)^2\int_{\R^2}\pk{\exists_{t\ge 0}:B_1(t)-t>x,B_2(t)-t>y}e^{(1-\rho^2)(x+y)}\td x\td y,$$
where $(B_1(t),B_2(t)),t \ge 0 $ is a 2-dimensional Gaussian process with $B_i$'s being standard Brownian motions with constant correlation $\tau$. Consequently, as $u\to\infty$
\bqny{
\psi_{2}(S,1,\vk 1u)\sim C_*u^{-2}e^{-\frac{u^2}{1+\tau} -\frac{c_*u}{2(1\CN{+}\tau)}},
}
where
\bqny{
C_*&=&\frac{e^{-\frac{c_*^{2}}{2(1-{\tau}^2)}}}{\sqrt{2\pi(1-{\tau}^2)}} \sum_{ i,j \in \{1\ldot d\}, \rho_{i,j}=\tau,~c_i+c_j=c_*}e^{\frac{c_ic_j}{1-\tau}}\\
&& \times \int_{\R^2}\pk{\exists_{t\ge 0}:B_1(t)-t>x,B_2(t)-t>y}e^{(1-{\tau}^2)(x+y)}\td x\td y \in (0,\IF).
}

{\bf Example 3} (Equi-correlated risk model): We consider the matrix $\Gamma$ such that $\Sigma=\Gamma
\Gamma^\top$ is an equi-correlated non-singular correlation matrix with off-diagonal entries equal to $\rho \in (-1/(d-1),1)$. Let $\vk a\inr^d$ with at least one positive component and assume for simplicity that its components are ordered, i.e., $a_1 \ge a_2\ge \cdots \ge  a_d$ and thus $a_1> 0$.
The inverse of $\Sigma$ equals
$$ \left[I_d - \vk 1 \vk 1^\top \frac{ \rho }{ 1+ \rho(d-1)} \right] \frac{1}{ 1-\rho},$$
where $I_d$ is the identity matrix. First we determine the index set $I$ corresponding to the unique solution of $\Pi_{\Sigma}(\vk a)$. We have for this case that $I$ with $m$ elements is unique and in view of \eqref{eq:IJ2}
\bqn{ \label{krue}
\vk \lambda_I= (\Sigma_{II})^{-1} \vk a_I= \frac 1 {1- \rho}\left[ \vk a_I - \rho \frac{\sum_{i\in I} a_i}{1+ \rho(m-1)} \vk 1_I\right] > \vk 0_I,
}
with $\vk 0 \inr^d$ the origin. From the above $m=\abs{I}=d$ if and only if
$$a_d> \rho \frac{\sum_{i=1}^d a_i}{1+ \rho(d-1)},$$
which holds in the particular case that all $a_i$'s are equal and positive.

When the above does not hold, the second condition on the index set $I$ given in \eqref{eq:IJ2} reads
$$ \Sigma_{JI}\Sigma_{II}^{-1} \vk a_I= \rho (\vk 1 \vk 1^\top )_{JI} \Sigma_{II}^{-1} \vk a_I \ge \vk a_{J}.$$
Next, suppose that $a_i=a>0, c_i=c\inr$ for all $i\le d$. In view of \eqref{lere} for any positive integer $k\le d$ and any $S\in [0,\Ee{1)}$ we have
\bqn{\psi_k(S, 1,\vk a u) &\sim & C\pk{\forall_{i \le k}: W_i(1) > ua+ c},
\quad u\to \IF,
\label{lere2}
}
where (set   $I=\{1 \ldot k\}$)
$$
C= \frac{d!}{k! (d-k)!£} \prod_{i\in I} \lambda_i \int_{\R^{k}} \pk{\exists_{t\ge 0}:\vk W_{I}(t)-t\vk a_{I}>\vk x_I}e^{ \vk \lambda_I^\top x_I} \td\vk x_I \in (0,\IF).
$$

Note that the case $\rho=0$ is treated in \cite{Krzyspp}[Prop. 3.6] and follows as a special case of this example.

\section{Proofs}

\subsection{\prooftheo{D}} Our proof below is based on the idea  of the proof of  \cite{KWW}[Thm 1.1],
where the case $\vk c$ has zero components, $k=d$ and $S=0$ has been considered.
Recall the definition of
\k1{sets $\vk E_\II(\vk u)$ and $E(\vk u)$ introduced in (\ref{def.E_I})}
for any non-empty $\II\subset\{1\ldot d\}$ such that $|\II|=k \le d$. With this notation we have
\bqny{
\psi_k(S, T,\vk u)=\pk{\exists{t\in[S,T]}:\vk W(t)-\vk c t\in\vk E(\vk u)} =
\pk{\tau_k(\vk u)\leq T},
}
where $\tau_k(\vk u)$ is the ruin time defined by
\bqny{
\tau_k(\vk u)=\inf\{ \CE{t \ge S}: \vk W(t)-\vk c t\in\vk E(\vk u)\}.
}
\kk{For the lower bound, we note that
\bqny{
\psi_k(S, T,\vk u)=\pk{\exists{t\in[S,T]}:\vk W(t)-\vk c t\in\vk E(\vk u)}\ge \pk{\vk W(T)-\vk c T\in\vk E(\vk u)}.
}
}
By the fact that Brownian motion has continuous sample paths
\bqn{\label{bF}
\vk W(\tau_k(\vk u))-\vk c \tau_k(\vk u)\in\partial\vk E(\vk u)
}
almost surely, where $\partial A $ stands for the topological boundary of the set $A \subset \R^d$.\\
Consequently, by the strong Markov property of the Brownian motion, we can write further
\bqny{
\lefteqn{\pk{\vk W(T) -\kk{\vk c}T \in\vk E(\vk u)}}\\
&= &\int_{0}^{T}\int_{\partial\vk E(\vk u)}\pk{\vk W(t)-\vk c t\in\td\vk x|\tau_k(\vk u)=t}\pk{\vk W(T)-\vk c\CN{T}\in\vk E(\vk u)|\vk W(t)-\vk c t=\vk x}\Ee{\pk{\tau_k(\vk u)\in\td t}}.
}
Crucial is that the boundary $\partial\vk E(\vk u)$ can be represented as the following union
\bqny{
\partial\vk E(\vk u)=\bigcup_{\substack{\II\subset\{1\ldot d\}\\|\II|=k}}\left(\partial\vk E_\II(\vk u)\cap\partial\vk E(\vk u)\right)=:\bigcup_{\substack{\II\subset\{1\ldot d\}\\|\II|=k}} F_\II(\kk{\vk u}).
}

For every $\vk x\in F_\II(\vk u)$ using the self-similarity of Brownian motion
\Ee{for all non-empty index sets $\mathcal{I} \subset \{1 \ldot d\}$ and all $t\in (\EZ{S},T)$}
\bqny{
\pk{\vk W(T)-\vk cT\in\vk E(\vk u)|\vk W(t)-\vk c t=\vk x}&\geq&\pk{\vk W(T)-\vk cT\in\vk E_\II(\vk u)|\vk W(t)-\vk c t=\vk x}\\
&=&\pk{\vk W_\II(T)-\vk c_\II T\geq \vk u_\II|\vk W(t)-\vk c t=\vk x}\\
&\geq&\pk{\vk W_\II(T-t)-\vk c_\II(T-t)\geq \vk 0}\\
&\geq&\pk{\vk W_\II(T-t)\geq \vk c_\II(T-t)}\\
&=&\pk{\vk W_\II(1)\geq \vk c_\II\sqrt{T-t}}\\
&\geq&\pk{\vk W_\II(1)\geq \tilde{\vk c}_\II\sqrt{T}}\\
&=&\pk{\vk W_\II(T)\geq \tilde{\vk c}_\II T}\\
&\geq&\min\limits_{\substack{\II\subset\{1\ldot d\}\\|\II|=k}}\pk{\vk W_\II(T)\geq \tilde{\vk c}_\II T},
}
where $\tilde{c}_i=\max(0,c_i)$, \Ee{hence
\k1{for all $\vk x \in \partial\vk E(\vk u)$}
$$  \pk{\vk W(\k1{T})-\vk c\k1{T}\in\vk E(\vk u)|\vk W(t)-\vk c t=\vk x}\geq
\min\limits_{\substack{\II\subset\{1\ldot d\}\\|\II|=k}}\pk{\vk W_\II(\k1{T})\geq \tilde{\vk c}_\II \k1{T}}.$$
}
 \k1{Consequently, \EZ{using further \eqref{bF} we obtain}}
\bqny{
& &\k1{\pk{\vk W({T}) -{\vk c}{T} \in\vk E(\vk u)}}\\
& &\qquad\k1{\ge}\min\limits_{\substack{\II\subset\{1\ldot d\}\\|\II|=k}}\pk{\vk W_\II(T)\geq \tilde{\vk c}_\II T}\int_{\EZ{S}}^{T}\int_{\partial\vk E(u)}\pk{\vk W(t)-\vk c t\in\td\vk x|\tau_k(\vk u)=t}\Ee{\pk{\tau_k(\vk u)\in\td t}}\\
& &\qquad=\min\limits_{\substack{\II\subset\{1\ldot d\}\\|\II|=k}}\pk{\vk W_\II(T)\geq \tilde{\vk c}_\II T}\psi_k(\EZ{S},T,\vk u)
}
establishing the proof.
\QED

\subsection{Proof of \netheo{theo1}} The results in this section hold under the assumption that
$\Sigma=\Gamma \Gamma^\top$ is positive definite, which \Ee{is equivalent with} our assumption that $\Gamma$ is non-singular.
The next lemma is a consequence of  \cite{MR3875603}[Lem 2].
\kk{We recall that
$\varphi$ denotes the probability density function of $ \kk{\Gamma}\vk B(1)$.}
\begin{lem}\label{proj1}
	For any $\vk a\in \R^d \setminus (-\IF, 0]^d$ we have for some positive constants $C_1, C_2$
	\bqny{
		\pk{\vk W(1) \CE{- \vk c }>\vk a u} \sim C_1 \pk{\forall_{i\in I}: W_i(1) \CE{- c_i }> a_i u}
		\sim C_2 u^{-\alpha}\varphi(\tilde{\vk a}u \CE{+ \vk c }), \quad u \to \IF,
	}
	where $\alpha$ is some \CE{integer} and $\tilde{\vk a}$ is the solution of quadratic programming problem $\Pi_{\Sigma}(\vk a), \Sigma=\Gamma \Gamma^\top $ and $I$ is the unique index set \Ee{that determines} the solution of $\Pi_{\Sigma}(\vk a)$. \Ee{ We agree in the following that
	if $\mathcal{I}$ is empty, then simply the term	$A_{\mathcal{I}}(t)$ should be deleted from the expressions below.}
\end{lem}

We state next  three lemmas  utilised in the case $T< \IF$.  Their proofs are displayed Section \ref{APP}. We recall that  $A_{\mathcal{I}}(t)$ is defined in \eqref{AI}.

\begin{lem}{\label{LL3:Lem}}
Let $\II,\JJ\subset\{1\ldot d\}$ be two index sets such that $\II\not=\JJ$ and $|\II|=|\JJ|=k \CE{\ge}1$.
If $\vk a_{\II\cup\JJ}$ has at least two positive components, then for any $s,t \in [0,1]$
\kk{there exists  some $\nu=\nu(s,t)>0$ such that as $u\to\infty$ }
\bqn{\label{LL3}
\pk{A_{\II}(t)\cap A_{\JJ}(s)}=\Ee{o\left(e^{-\nu u^2}\right)}\sum_{\substack{\II^*
\subset\{1\ldot d\}\\|\II^*|=k}}\pk{A_{\II^*}(1)}, 
}
\Ee{and}
\bqn{\label{LL3.1}
\pk{A_{\II\setminus\JJ}(t), A_{\kdd{\JJ\setminus\II}}(s), A_{\II\cap\JJ}(\min(t,s))}=\Ee{o\left(e^{-\nu u^2}\right)}\sum_{\substack{\II^*\subset\{1\ldot d\}\\|\II^*|=k}}\pk{A_{\II^*}(1)}.
}
\end{lem}

\begin{lem}{\label{LL4.1:Lem}}  Let $S > 0$, $k\le d$ be a positive integer and let $\vk a \inr^d$ be given. If $\II,\JJ\subset\{1\ldot d\}$ are two \Ee{different} index sets with $k\ge 1$ elements such that $\vk a_{\II\cup\JJ}$ has at least one positive component, then there exist $s_1,s_2\in [S,1]$ and some positive constant $\tau$ such that as $u\to\infty$
\bqn{
\pk{\exists s,t\in[S,1]: A_{\II}(s)\cap A_{\JJ}(t)}=\Ee{o\left(e^{\tau u}\right)} \pk{A_{\II\setminus\JJ}(s_1)\cap A_{\JJ\setminus\II}(s_2)\cap A_{\II\cap\JJ}(\min(s_1,s_2))}.\label{LL4.1}
}
\end{lem}

\underline{\bf Case $T< \IF$.} According to \netheo{D} and  \nelem{proj1} it is enough to show the proof for $S\in \cE{(0,T)}$. In view of the self-similarity of Brownian motion we assume for simplicity $T=1$.
Recall that in our notation $\Sigma= \Gamma  \Gamma ^\top$
is the covariance matrix of $\vk W(1)$ which is non-singular  and we denote its pdf by $\varphi$. In view of \eqref{LL3}, \eqref{LL3.1} and \eqref{LL4.1}  for all $S\in (0,1)$ \kdd{there exists some $\nu>0$ such that} \kk{as $u\to\infty$}
\bqny{
\sum_{\substack{\II,\JJ \subset\{1\ldot d\}\\| \II|=|\JJ|=k, \II\not= \JJ}}\pk{\exists s,t\in[S,1]: A_\II(s)\cap A_\JJ(t)}
=\kk{o\left( e^{- \nu u^2} \right)}\sum_{\substack{\II\subset\{1\ldot d\}\\|\II|=k}}\pk{A_\II(1)}.
}

Note that we may  \Ee{utilise \eqref{LL3.1}  and \eqref{LL4.1}}
for sets $\II$ and $\JJ$ of length $k$, because of the assumption that $\vk a$ has no more than $k-1$ non-positive components. Hence any vector $\vk a_{\II}$ has at least one positive component.

Further, by \kk{Theorem} \ref{D} and the inclusion-exclusion formula we have that for
\kk{some $K>0$ and}  all $u$ sufficiently large
$$ \psi_k(S, 1,\vk u) \le \kk{K}\sum_{\substack{I\subset\{1\ldot d\}\\|I|=k}}\pk{A_\II(1)}.
$$
Hence the claim follows from \eqref{inclusoin1} and \eqref{inclusoin2}.

\underline{\bf Case $T=\IF$}. \COM{In this case, according to \netheo{D} and  \nelem{proj1},
it is enough to show the proof
 for the case $S=0$.}
Using the self-similarity of Brownian motion we have 
\bqny{
\pk{\exists t>0: A_{\II}(t)}=\pk{\exists t>0:\vk W_\II(ut)\ge(\vk a+\vk c t)_{\kk{\II}} u}
&=&\pk{\exists t>0:\vk W_\II(t)\ge (\vk a+\vk c t)_{\kk{\II}}\sqrt{u}}\\
&=&\pk{\exists t>0:A^*_\II(t)},
}
\Ee{where
\bqn{\label{AY}
	A^*_\II(t)= \{   \vk W_\II(t)\ge (\vk a+\vk c t)_{\kk{\II}}\sqrt{u} \} .}
}
For $t>0$ define 
\bqn{
r_\II(t)=\min_{\vk x\ge\vk a_{\II}+\vk c_{\II} t}\frac{1}{t}\vk x^\top\Sigma^{-1}_{\II\II}\vk x,\label{rfunct}
\ \ \Sigma_{\II\II}=Var(\vk W_{\II}(1)), \ \ \Sigma^{-1}_{\II\II}= (\Sigma_{\II\II})^{-1}.
}
Since  $\lim_{t\downarrow 0}r_{\kk{\II}}(t)=\infty$ we set below    $r_{\kk{\II}}(0)=\infty$.\\

In view of \kk{Lemma} \ref{proj1} we have as   $u\to \IF$
\bqny{
\pk{A_\II^*(t)}\sim C_1 u^{-\alpha/2}\varphi_{\II,t}((\widetilde{\vk a_{\II}+\vk c_{\II} t})\sqrt{u})=
 C_2 u^{\kk{-\alpha/2}} e^{-\frac{r_\II(t) u}{2}},
}
where $\widetilde{\vk a_{\II}+\vk c_{\II} t}$ is the solution of quadratic programming problem
$\Pi_{t\Sigma_{\II\II}}(\vk a_{\II}+\vk c_{\II} t)$
and $\varphi_{\II,t}(\vk x)$ is the pdf  of $\vk W_{\II}(t)$, $\alpha$ is some integer
and \kdd{$C_1,C_2$ are positive constant that  do not depend on $u$}. 
For \kk{notational simplicity}  we shall omit below the subscript $\II$.

The \kdd{case  $T=\infty$ is} established utilising  the following two lemmas,
whose proofs are \kdd{displayed} in Section \ref{APP}.

\begin{lem}\label{diff}
Let $k\le d$ be a positive integer and let $\vk a,\vk c\in\R^d$. Consider two different sets $\II,\JJ\subset\{1\ldots d\}$ of cardinality $k$. If both $\vk a_\II+\vk c_\II t$ and $\vk a_\JJ+\vk c_\JJ t$ have at least one positive component for all $t>0$ and both $\vk c_{\II}$ and $\vk c_\JJ$ also have at least one positive component, then in case $\hat{t}_\II:=arg\min\limits_{t>0}~r_{\II}(t)\not=\hat{t}_\JJ:=arg\min\limits_{t>0}~r_{\JJ}(t)$,
\bqny{
\pk{\exists s,t>0:~A^*_{\II}(t)\kk{\cap}  A^*_{\JJ}(s)}=o(\pk{A^*_\II(\hat{t}_{\II})}+\pk{A^*_\JJ(\hat{t}_\JJ)}), \quad u\to \IF.
}
\end{lem}

\begin{lem}\label{same}
\Ee{Under the settings of \nelem{diff}, if} $\vk a+\vk c t$ has no more than $k-1$ non-positive component for all $t>0$ and $\CNN{\vk c}$ has no more than $k-1$ non-positive components, then in case $\hat{t}_\II:=arg\min\limits_{t>0}~r_{\II}(t)=\hat{t}_\JJ:=arg\min\limits_{t>0}~r_{\JJ}(t)$
\bqny{
\pk{\exists s,t>0:~A^*_{\II}(t)\kk{\cap} A^*_{\JJ}(s)}=o\left(\sum_{\substack{\KK\subset\{1\ldots d\}\\|\KK|=k}}\pk{A^*_\KK(\hat{t}_{\KK})}\right), \quad u\to \IF.
}
\end{lem}

Combining the above two lemmas we have that for any two index sets $\II,\JJ\subset\{1\ldot d\}$ of cardinality $k$, there is some index set $\KK\subset\{1\ldot d\}$ such that as $u\to \IF$
\bqny{
\pk{\exists s,t>0: A_{\II}^*(s)\kk{\cap} A_{\JJ}^*(t)}=o\left(\pk{\exists t>0:A_{\KK}^*(t)}\right),
}
which is \Ee{equivalent with}
\bqny{
\pk{\exists s,t>0: A_{\II}(s)\kk{\cap} A_{\JJ}(t)}=o\left(\pk{\exists t>0:A_{\KK}(t)}\right).
}
The proof follows now by \eqref{inclusoin1} and \eqref{inclusoin2}.
 \QED

\subsection{Proof of \netheo{theo2}}

Below we set
$$ \delta(u,\CN{\Lambda}):=1-\CN{\Lambda}u^{-2}$$
and denote by $\tilde{ \vk a}$ the unique solution of the quadratic programming problem $\Pi_\Sigma(\vk a)$.

We denote below by $I$ the index set that determines the unique solution of $\Pi_{\Sigma} (\vk a)$,
where $\vk a \inr^d$ has at least one positive component (see Lemma \ref{prop1}). If $J= \{1 \ldot d \} \setminus I$ is non-empty, then we set below $U=\{j\in J:\tilde{ a}_j= a_j\}$. The number of elements $|I|$ of $I$ is denoted by $m$, which is a positive integer.

The next \CN{lemma} \EZ{is} proved in  Section \ref{APP}.

\BEL\label{inimp:Lem}
For any $\CN{\Lambda}>0$, $\vk a\in \R^d \setminus (-\infty,0]^d, \Ee{\vk c \inr^d}$ and all sufficiently large u there exist $C>0$ such that
\bqn{\label{inimp}
m(u,\CN{\Lambda}):=\pk{\exists_{t\in[0,\delta(u,\CN{\Lambda})]}:\vk W(t)-t\vk c> u\vk a}\le e^{-\CN{\Lambda}/C}\frac{\pk{\vk W(1)\ge \vk a u+\vk c}}{\pk{\vk W(1)>\max(\vk c,0)}}
}
\Ee{and further}
\bqn{\label{A}
M(u,\CN{\Lambda}):=\pk{\exists_{t\in[\delta(u,\CN{\Lambda}),1]}:\vk W(t)-t\vk c >u\vk a}\sim
C(\vk c)K([0,\CN{\Lambda}])u^{-m}\varphi(u\tilde{\vk a}+\vk c),
}
where $C(\vk c)= \pk{ \vk W_U(1)> \vk c_U \lvert
\vk W_I(1)> \vk c_I } $ and for $\vk \lambda= \Sigma^{-1} \tilde{\vk a}$
\bqny{
E([\CN{\Lambda_1},\CN{\Lambda_2}])=\int_{\R^{m}} \pk{\exists_{t\in[\CN{\Lambda_1},\CN{\Lambda_2}]}:\vk W_{I}(t)-t\vk a_{I}>\vk x_I}e^{ \vk\lambda_I ^\top \vk x_I}\td\vk x_I \in (0,\IF)
}
for all   constants $\CN{\Lambda_1} < \CN{\Lambda_2}$. We set $C(\vk c)$ equal 1 if $U$ defined in Remark \ref{rem1}  is empty. Further we have
\bqn{ \label{T}
\lim\limits_{\CN{\Lambda}\to\infty}E([0,\CN{\Lambda}])=\int_{\R^{m}}\pk{\exists_{t\geq 0}:\vk W_{I}(t)-t\vk a_{I}>\vk x_I}e^{\vk \lambda_I^\top \vk x_I}\td\vk x_I\in(0,\infty).
}

\EEL

First note that for all $\CN{\Lambda},u$ positive
\bqny{
M(u,\CN{\Lambda})\leq\pk{\exists_{t\in[0,1]}:\vk W(t)-t\vk c>u\vk a}\leq M(u,\CN{\Lambda})+m(u,\CN{\Lambda}).
}
In view of Lemmas \ref{inimp:Lem} and \ref{proj1}
\bqny{
\lim_{\CN{\Lambda}\to\infty}\lim_{u\to\infty}\frac{m(u,\CN{\Lambda})}{M(u,\CN{\Lambda})}=0,
}
hence
\bqny{
\lim_{\CN{\Lambda}\to\infty}\lim_{u\to\infty}\frac{\pk{\exists_{t\in[0,1]}:\vk W(t)-t\vk c>u\vk a}}{M(u,\CN{\Lambda})}=1
}
and thus the proof follows applying \eqref{A}. \QED
\COM{
\bqny{
\lim_{T\to\infty}\lim_{u\to\infty}\frac{\pk{\exists_{t\in[0,1]}:\vk W(t)-t\vk c>u\vk a}}{M(u,T)}&=&\lim_{T\to\infty}\lim_{u\to\infty}\frac{\pk{\exists_{t\in[0,1]}:\vk W(t)-t\vk c>u\vk a}}{ Q K([0,T])F(u)}\frac{Q K([0,T])F(u)}{M(u,T)}\\
&=&\lim_{T\to\infty}\lim_{u\to\infty}\frac{\pk{\exists_{t\in[0,1]}:\vk W(t)-t\vk c>u\vk a}}{ QK([0,T])F(u)}.
}
Finally, using \eqref{T} we obtain
\bqny{
\lim_{T\to\infty}\lim_{u\to\infty}\frac{\pk{\exists_{t\in[0,1]}:\vk W(t)-t\vk c>u\vk a}}{ Q K([0,T])F(u)}=\lim_{u\to\infty}\frac{\pk{\exists_{t\in[0,1]}:\vk W(t)-t\vk c>u\vk a}}{ Q K([0,T]) F(u)}=1
}
establishing the proof. \QED
}

\subsection{Proof of Eq. \eqref{time2}}
The proof is similar to that of  \cite{DHL17}[Thm 2.5] and therefore we highlight only the main steps that lead to \eqref{time2}.

Let $T>S\ge 0$. According to the definition of $\tau(u)$ and the self-similarity of
Brownian motion
\bqny{
\frac{\tau(u)}{T}=\inf\{t\ge 0:\vk W(Tt)-tT\vk c>\vk a u\}=\inf\{t\ge 0:\vk W(t)-t\sqrt{T}\vk c>\vk a u/\sqrt{T}\}.
}
Thus, without loss of generality in the rest of the proof we suppose that $T=1>S\ge0$.

We note that
\bqny{
\lefteqn{
\pk{u^2(1-\tau(u))\geq x| \tau(u)\in [S,1]}=\frac{\pk{u^2(1-\tau(u))\geq x, \tau(u)\in [S,1]}}{\pk{\tau(u)\in [S,1]}}}\\
&=& \frac{\pk{u^2(1-\tau(u))\geq x, \tau(u)\leq 1}}{\pk{\tau(u)\in [S,1]}}
-\frac{\pk{u^2(1-\tau(u))\geq x, \tau(u)\leq S}}{\pk{\tau(u)\in [S,1]}}\\
&=&P_1(u)-P_2(u).
}
Next, for $\tilde{x}(u)=1-\frac{x}{u^2}$
\begin{eqnarray*}
P_1(u)&=&\frac{\pk{\tau(u)\leq \tilde{x}(u)}}{\pk{\tau(u)\in [S,1]}}\sim
\frac{\pk{\exists_{t\in[0,\tilde{x}(u)]}:\vk W(t)-\vk c t>u\vk a}}
     {\pk{\exists_{t\in[0,1]}:\vk W(t)-\vk c t>u\vk a}}\\
&=&
\frac{\pk{\exists_{t\in[0,1]}:\vk W(t)-(\vk c\sqrt{\tilde{x}(u)}) t>\frac{u}{\sqrt{\tilde{x}(u)}}\vk a}}
     {\pk{\exists_{t\in[0,1]}:\vk W(t)-\vk c t>u\vk a}}, \ u \to \IF.
\end{eqnarray*}
Hence Theorem \ref{theo2} and
\bqny{
\varphi\left(\frac{u}{\sqrt{\tilde{x}(u)}}\tilde{\vk a}+(\vk c\sqrt{\tilde{x}(u)})\right)=\varphi(u\tilde{\vk a}+\vk c)e^{-\frac{1}{2}\left(\frac{1}{\tilde{x}(u)}-1\right)u^2\tilde{\vk a}^\top\Sigma^{-1}\tilde{\vk a}}e^{-\frac{1}{2}(\tilde{x}(u)-1)\vk c^\top\Sigma^{-1}\vk c}
}
and
\bqny{
\limit{u} e^{-\frac{1}{2}\left(\frac{1}{\tilde{x}(u)}-1\right)u^2\tilde{\vk a}^\top\Sigma^{-1}\tilde{\vk a}}&=& e^{-x\frac{\tilde{\vk a}^\top\Sigma^{-1}\tilde{\vk a}}{2}},\quad
\limit{u}e^{-\frac{1}{2}(\tilde{x}(u)-1)\vk c^\top\Sigma^{-1}\vk c}= 1
}
we obtain
\begin{eqnarray}
\lim_{u\to\infty} P_1(u)= e^{-x\frac{\tilde{\vk a}^\top\Sigma^{-1}\tilde{\vk a}}{2}}.\label{P1}
\end{eqnarray}
Moreover, following the same reasoning as above
\begin{eqnarray}
P_2(u)=\frac{\pk{ \tau(u)\leq S}}{\pk{\tau(u)\in [S,1]}}
&\sim&
\frac{\pk{ \tau(u)\leq S}}{\pk{\tau(u)\leq 1}}\to 0
\label{P2}
\end{eqnarray}
as $u\to\infty$.
Thus, combination of (\ref{P1}) with (\ref{P2}) leads to
\bqn{\label{time0}
\limit{u} \pk{u^2(1-\tau(u))\geq x|\tau(u) \in [S,1]} = e^{-x\frac{\tilde{\vk a}^\top\Sigma^{-1}\tilde{\vk a}}{2}}.\nonumber
}
\QED

\section{Appendix}
\label{APP}
\begin{lem}{\label{LL1}}
\Ee{If for $\vk a\in (\R\cup\{-\infty\})^d$ and $\II\subset\{1\ldot d\}$} such that $\vk a_\II$ has at least two positive components and $\Gamma$ is non-singular, then for all $t>0$
\bqny{
\pk{A_\II(t)}=\Ee{o\left(e^{-\nu u^2}\right)}\sum_{i\in \II}\pk{A_{ \II\setminus\{i\}}(t)}, \quad u\to \IF,
}
where $\nu=\nu(t,\II)>0$ does not depend on $u$.
\end{lem}
\begin{remark}\label{reduct}
Lemma \ref{LL1} implies that for any vector $\vk a\in (\R\cup\{-\infty\})^d$ and for any $d$-dimensional Gaussian random vector $\vk W$, if $\vk a$ has at least two positive components, there exists some positive constant $\eta$ and $i\in\{1\ldots d\}$ such that as $u\to \IF$
\bqny{
\pk{\vk W>\vk a u}=o(e^{-\eta u^2})\pk{\vk W_{K}>\vk a_{K}u}, \quad K= \{1\ldot d\}\setminus \{i\} .
}
\end{remark}

\prooflem{LL1} For \kk{notational simplicity}  we shall assume that $\II= \{1\ldot d \}$ and set $K_i=\II\setminus \{i\} $. By the assumption for all $i\in\II$ the vector $\vk a_{K_i}$ has at least one positive component and $\Sigma=\Gamma\Gamma^\top$ is positive definite. In view of \nelem{proj1} for any fixed $t>0$
and some $C_1,C_2$ two positive constants we have
\bqny{
\pk{A_\II(t)}\sim C_1 u^{\alpha_1}\varphi_t(\tilde{\vk a}u+\vk c), \quad
\pk{A_{ K_i}(t)} \sim C_2 u^{\alpha_2}\varphi_t(\bar{\vk a_i}u +\vk c), \quad u\to \IF,
}
where $\varphi_t$ is the pdf of   $\vk W(t)$ with covariance matrix $\Sigma(t)=t\Sigma$ and
\bqny{
\tilde{\vk a}=arg\min_{\vk x \ge \vk a} \vk x^\top \Sigma^{-1}(t)\vk x, \quad
\bar{\vk a_i}=arg\min_{\vk x\in S_i} \vk x^\top\Sigma^{-1}(t)\vk x,
}
with $$S_i=\{\vk x\in\R^{d}:~\forall j \in K_i : x_j\geq a_j\}.$$
Since $\{\x\inr^d: \x \ge \vk a \}\subset S_i$, then clearly
$$\tilde{\vk a}^\top \Sigma^{-1}(t) \tilde{\vk a} \geq\bar{\vk a_i}^\top \Sigma^{-1}(t)\bar{\vk a_i}$$
 for any $i\le d$.
Next, if we have strict inequality for some $i\le d$, i.e., $\tilde{\vk a}^\top\Sigma^{-1}(t) \tilde{\vk a}>\bar{\vk a_i}^\top \Sigma^{-1}(t)\bar{\vk a_i}$, then it follows that
\bqny{
\pk{A_\II(t)}\sim C u^{\alpha_1}\varphi_t(\tilde{\vk a}u +\vk c)=o\left(e^{-\nu u^2}\pk{A_{K_i}(t) }\right), \quad u\to \IF
}
for $\nu=\frac{1}{2}\left(\tilde{\vk a}^\top\Sigma^{-1}(t) \tilde{\vk a}-\bar{\vk a_i}^\top \Sigma^{-1}(t)\bar{\vk a_i}\right) >0$, hence the claim follows.\\
\COM{Let us consider now the extreme case that for all $i\le d$ we have $\tilde{\vk a}^\top\Sigma^{-1}\tilde{\vk a}=\bar{\vk a_i} ^\top \Sigma^{-1}\bar{\vk a_i}$.
Since from the following well-known identity valid for any two non-overlaping index set $A,B, A\cup B= \{ 1 \ldot d\}$  and any $\vk x\inr$  (using Schur compliments and noting  that $(\Sigma^{-1})_{BB} $ is positive definite)
\bqny{
	\x^\top  \Sigma^{-1} \x =\x_A^\top (\Sigma_{AA})^{-1} \x_A  +
	(\x_B - \Sigma_{BA}(\Sigma_{AA})^{-1}  \x_A) ^\top (\Sigma^{-1})_{BB}  (\x_B - \Sigma_{BA}(\Sigma_{AA})^{-1}  \x_A)
}
it follows that (set $K=K_i$)
$$ \bar{\vk a_i}=arg\min_{\vk x_{K} \ge \vk a_{K} } \vk x^\top_K\CE{(\Sigma_{KK})^{-1}}\vk x_K.$$

In view of \nelem{proj1} the asymptotics of $\pk{ A_\II(t)}$ is determined by the unique index set $I$ that
determines $\bar{\vk a}$. Since $\bar{\vk a_i}$ is equal to $\bar{\vk a}$ for all $i\le d$, then the
asymptotics of $\pk{A_{K_i}(t)}$ is up to some constants the same as that of  $\pk{ A_\II(t)}$ as $u\to \IF$.
 This implies by \nelem{proj1} that the corresponding index sets $I_i$'s of the quadratic programming problems must agree for all
 $i\le d$ with $I$, which is a contradiction.\\}
Let us consider now the extreme case that for all $i\le d$ we have $\tilde{\vk a}^\top\Sigma^{-1}\tilde{\vk a}=\bar{\vk a_i} ^\top \Sigma^{-1}\bar{\vk a_i}$. As we know that each $\bar{\vk a}_i$ is unique, then $\bar{\vk a}_i=\tilde{\vk a}$ for all $i\in\II$.
Consider set
\bqny{
E=\{\vk x\inr^d:\vk x^\top\Sigma^{-1}(t)\vk x\leq\tilde{\vk a}^\top \Sigma^{-1}(t) \tilde{\vk a}\}.
}
\kdd{Since} $\Sigma(t)$ is positive definite, $E$ is a full dimensional ellipsoid in $\R^d$. By the definition,
$E\cap S_i=\{\tilde{\vk a}\}$. Define the following lines in $\R^d$
\bqny{
l_i=\{\vk x\inr^d:\forall i\in K_i, x_i=\tilde{a}_i\}
}
\kdd{and observe that since $l_i\in S_i$, then} $l_i\cap E=\{\tilde{\vk a}\}$, and they are linearly independent.
Since $E$ is smooth, there can not be more than $d-1$ linearly independent tangent lines at the point $\tilde{\vk a}$,
\kdd{which leads to} a contradiction.

 \QED

\prooflem{LL3:Lem}
First note that since $\II\not= \JJ$, then $|\II\cup \JJ|\geq k+1$. Consequently, we can find some \Ee{index} set $\KK$ such that
$$|\KK|=k+1, \quad \KK\subset \II\cup \JJ$$
and further $ \vk a_ {\KK}$ has at least two positive components. Applying \nelem{LL1} for any $t\in [0,1]$ and some $\nu>0$
\bqny{
\pk{A_{\KK}(t)}= o\left(e^{-\nu u^2}\right)\sum_{j\in \KK}\pk{A_{\KK\setminus \{j\}}(t)}, \quad u\to \IF.
}
If $s=t$, then applying \nelem{proj1}
\bqny{
0\leq\pk{A_{\II}(t)\cap A_{\JJ}(t)}=\pk{A_{\II\cup \JJ}(t)}\leq\pk{A_{\KK}(t)}= o\left(e^{- \nu u^2} \right)\sum_{\substack{\II^*\subset\{1\ldot d\}\\ \II^*|=k}}\pk{A_\II^*(t)}.
}
Next, if $s<1$, then applying \nelem{proj1} we obtain
\bqny{
0\leq\pk{A_{\II}(t)\cap A_{\JJ}(s)}\leq\pk{A_{\JJ}(s)}=o\left(e^{-\nu u^2}\pk{A_{\JJ}(1)}\right)=o\left(e^{-\nu u^2}\right)\sum_{\substack{\II^*\subset\{1\ldot d\}\\|\II^*|=k}}\pk{A_{\II^*}(1)}.
}
A similar asymptotic bound follows for $t<1$, whereas if $s=t=1$, the \Ee{first} claim follows directly from the case $s=t$ discussed above.

\Ee{We show next \eqref{LL3.1}.}
If $s<t$, then $s<1$ and applying \nelem{proj1} we obtain
\bqny{
0 &\le& \pk{A_{\II\setminus\JJ}(t), A_{\JJ\setminus\II}(s), A_{\II\cap\JJ}(\min(t,s))}\\
&\leq &\pk{A_{\JJ}(s)}=o\left(e^{-\nu u^2}\pk{A_{\JJ}(1)}\right) \\
&=&o\left(e^{-\nu u^2} \right)\sum_{\substack{\KK\subset\{1\ldot d\}\\|\KK|=k}}\pk{A_{\KK}(1)}.
}
A similar asymptotic bound follows for $t<s$ or $s=t\le 1$ by applying \eqref{LL3} establishing the proof.
\QED

\prooflem{LL4.1:Lem}
Define for $s,t\in[S,1]$ the Gaussian random vector
$$\mathcal{W}(s,t)=(\vk W_{\II\setminus\JJ}(s)^\top, \vk W_{\JJ\setminus\II}(t)^\top, \vk W_{\II\cap\JJ}(\min(s,t))^\top)^\top,$$
with covariance matrix $D(s,t)$. We show first that this matrix is positive definite. For this we assume that $s\leq t$. As $D(s,t)$ is some covariance matrix, we know that it is non-negative definite. Choose some vector $\vk v\in\R^d$. It is sufficient to show that if $\vk v^\top D(s,t)\vk v=\vk 0$, then $\vk v=\vk 0$ (here $\vk 0=(0, \ldots,0)^\top\in\R^d$). Note that
$$
\vk v^\top D(s,t)\vk v=Var(\sprod{\mathcal{W}(s,t),\vk v})=Var(\sprod{\vk W(s),\vk v}+\sprod{\vk W_{\JJ\setminus\II}(t)-\vk W_{\JJ\setminus\II}(s),\vk v_{\JJ\setminus\II}}).
$$
Using that $\vk W(t)$ has independent increments, this variance is equal to the sum of the variances. Hence, both of them should be equal to zero. In particular it means that $Var(\sprod{\vk W(s),\vk v})=0$. Hence, as $s\geqslant S>0$, we have that $\vk v=\vk 0$. Thus, $D(s,t)$ is positive definite and $D^{-1}(s,t)$ exists. \\
Set further
$${\mathfrak{a}=(\vk a_{\II\setminus\JJ}^\top,\vk a_{\JJ\setminus\II}^\top, \vk a_{\II\cap\JJ}^\top)^\top}, \quad {\mathfrak{c}(s,t)=(s\vk c_{\II\setminus\JJ}^\top, t\vk c_{\JJ\setminus\II}^\top, \min(s,t)\vk c_{\II\cap\JJ}^\top)^\top}.$$ With this notation we have
\bqny{
\pk{\exists s,t\in[S,1]:A_{\II}(s)\cap A_{\JJ}(t)}\leq\pk{\exists s,t\in[S,1]:\mathcal{W}(s,t)-\mathfrak{c}(s,t)\geq\mathfrak{a}u}.
}
Let $\tilde{\mathfrak{a}}(s,t) = arg\min_{\kk{\vk x}\geq\mathfrak{a}} \vk x^\top D^{-1}(s,t) \vk x$ be the unique solution of $\Pi_{D(s,t)}(\mathfrak{a})$ and let further $
\mathfrak{w}(s,t) =D ^{-1}(s,t)\tilde{\mathfrak{a}}(s,t)$
be the solution of the dual problem. We denote by $I(s,t)$ the index set related to the quadratic programming problem
$\Pi_{D(s,t)}(\mathfrak{a})$. Then $ \mathfrak{w}(s,t)$ has non-negative components and according to \nelem{prop1} since  both
$s,t\geq S>0$ we have
$$\mathfrak{a}^\top\mathfrak{w}(s,t)=\tilde{\mathfrak{a}}^\top(s,t)\mathfrak{w}(s,t) =\tilde{\mathfrak{a}}^\top(s,t) D^{-1}(s,t)
\tilde{\mathfrak{a}}(s,t)>0.$$
 Consequently, we have
\bqny{
\pk{\exists s,t\in[S,1]:\mathcal{W}(s,t)-\mathfrak{c}(s,t)\geq\mathfrak{a}u}&\leq&\pk{\exists s,t\in[S,1]: \mathfrak{w}^\top(s,t) \left(\mathcal{W}(s,t)-\mathfrak{c}(s,t)\right)\geq u \mathfrak{w}^\top(s,t) \tilde{\mathfrak{a}}(s,t)}\\
&=&\pk{\exists s,t\in[S,1]:\frac{\mathfrak{w}^\top(s,t) \left(\mathcal{W}(s,t)-\mathfrak{c}(s,t)\right)}{\mathfrak{w}^\top(s,t) \tilde{\mathfrak{a}}(s,t)}\geq u}\\
&\leq&\pk{\exists s,t\in[S,1]:\frac{\mathfrak{w}^\top(s,t) \mathcal{W}(s,t)}{\mathfrak{w}^\top(s,t) \tilde{\mathfrak{a}}(s,t)}\geq u+\mathfrak{C}}
}
for any positive $u$, where $\mathfrak{C}=\min_{s,t\in[S,1]}\frac{\mathfrak{w}^\top(s,t)\mathfrak{c}(s,t)}{\mathfrak{w}^\top(s,t)\tilde{\mathfrak{a}}(s,t)}$. Moreover,
for some $s_1,s_2\in[S,1]$
\bqny{
\sigma^{2}=\sup_{s,t\in[S,1]}\E{\left(\frac{{\mathfrak{w}^\top(s,t)\mathcal{W}(s,t)}}{{\mathfrak{w}^\top(s,t) \tilde{\mathfrak{a}}(s,t)}}\right)^2}
=\sup_{s,t\in[S,1]}\frac{1}{\tilde{\mathfrak{a}}^\top(s,t) D^{-1}(s,t)\tilde{\mathfrak{a}}(s,t)}
=\frac{1}{\tilde{\mathfrak{a}}^\top(s_1,s_2)D^{-1}(s_1,s_2)\tilde{\mathfrak{a}}(s_1,s_2)}
}
since $[S,1]^{2}$ is compact. In order to utilize Piterbarg inequality, see e.g., \cite{Pit96}[Thm 8.1],
we need to show that for some positive constant $G$
$$\E{\left(\frac{\mathfrak{w}^\top(s_1,t_1) \mathcal{W}(s_1,t_1)}{\mathfrak{w}^\top(s_1,t_1) \tilde{\mathfrak{a}}(s,t)}-\frac{\mathfrak{w}^\top(s_2,t_2) \mathcal{W}(s_2,t_2)}{\mathfrak{w}^\top(s_2,t_2) \tilde{\mathfrak{a}}(s,t)}\right)^2}\leq G
[\abs{s_1-s_2}+\abs{t_1-t_2}],\qquad s_1,s_2,t_1,t_2\in[S,1].$$
It is enough to show that for some positive constant $G_1,G_2$ and all $s_1,s_2,t_1,t_2\in[S,1]$
\Ee{(set $\vk 1:=(1\ldot 1)^\top\in\R^d$)} $$\E{\left(\mathcal{W}(s_1,t_1)-\mathcal{W}(s_2,t_2)\right)^2}\leq G_1\vk 1 [\abs{s_1-s_2}+\abs{t_1-t_2}]
$$ and
\bqn{
\E{\left(\frac{\mathfrak{w}(s_1,t_1)}{\mathfrak{w}^\top(s_1,t_1) \mathfrak{a}}-\frac{\mathfrak{w}(s_2,t_2)}{\mathfrak{w}^\top(s_2,t_2) \mathfrak{a}}\right)^2}\leq G_2\vk 1[\abs{s_1-s_2}+\abs{t_1-t_2}].
\label{holder}
}
The first inequality is clear from the definition of $\mathcal{W}$. Consider the second one. For any real \Ee{square}
matrix $M$ below we write
$\lvert\lvert M\rvert\rvert_{\sup}$ for the maximal absolute value of its elements.
Note that for some positive constant $G^*$
$$\left\|D(s_1,t_1)-D(s_2,t_2)\right\|_{\sup}\leq G^*
[\abs{s_1-s_2}+\abs{t_1-t_2}].$$
It implies that according to \cite{Hag79}[Thm 3.1] for some constant $G^\prime$
$$\|\tilde{\mathfrak{a}}(s_1,t_1)-\tilde{\mathfrak{a}}(s_2,t_2)\|_{\sup}\leq G^\prime [\abs{s_1-s_2}+\abs{t_1+t_2}]$$
and also $$\bigl|\abs{D(s_1,t_1)}-\abs{D(s_2,t_2)}\bigr|\leq\abs{d}!dG^*\sup_{s,t\in[S,1]}\|D(s,t)\|_{\sup}^{d-1}(\abs{s_1-s_2}+\abs{t_1-t_2}).$$
Thus (write $C(s,t)$ for \CNN{the adjugate}  matrix $D(s,t)$, i.e., it is equal to the transpose of the cofactor matrix of $D(s,t)$)
$$\left\|C(s_1,t_1)-C(s_2,t_2)\right\|_{\sup}\leq(d-1)!(d-1)G^*\sup_{s,t\in[S,1]}\|D(s,t)\|_{\sup}^{d-2}(\abs{s_1-s_2}+\abs{t_1-t_2}).$$
Since  $\abs{D(s,t)}$ is strictly positive for $s,t\in[S,1]$ (i.e., $\abs{D(s,t)}>\delta$ for some positive $\delta$ and all $s,t\in[S,1]$), then
\bqn{
\lefteqn{\left\|\mathfrak{w}(s_1,t_1)-\mathfrak{w}(s_2,t_2)\right\|_{\sup}}\\
& =&\left\|\abs{D(s_1,t_2)}^{-1}C(s_1,t_1)\tilde{\mathfrak{a}}(s_1,t_1)-\abs{D(s_2,t_2)}^{-1}C(s_2,t_2)\tilde{\mathfrak{a}}(s_2,t_2)\right\|_{\sup}\notag \\
& \leq &d\abs{\abs{D(s_1,t_1)}^{-1}-\abs{D(s_2,t_2)}^{-1}}\left\|C(s_1,t_1)\right\|_{\sup}\left\|\tilde{\mathfrak{a}}(s_1,t_1)\right\|_{\sup}\notag \\
& &+d\abs{\abs{D(s_2,t_2)}^{-1}}\left\|C(s_1,t_1)-C(s_2,t_2)\right\|_{\sup}\left\|\tilde{\mathfrak{a}}(s_1,t_1)\right\|_{\sup}\notag \\
& &+d\abs{\abs{D(s_2,t_2)}^{-1}}\left\|C(s_2,t_2)\right\|_{\sup}\left\|\tilde{\mathfrak{a}}(s_1,t_1)-\mathfrak{a}(s_2,t_2)\right\|_{\sup}\notag \\
&\leq& \frac{d}{\delta^2}\abs{\abs{D(s_1,t_1)}-\abs{D(s_2,t_2)}}\sup_{s,t\in[S,1]}\left\|C(s,t)\right\|_{\sup}\sup_{s,t\in[S,1]}\left\|\tilde{\mathfrak{a}}(s,t)\right\|_{\sup}\notag \\
& &+\frac{d}{\delta}\left\|C(s_1,t_1)-C(s_2,t_2)\right\|_{\sup}\sup_{s,t\in[S,1]}\left\|\tilde{\mathfrak{a}}(s,t)\right\|_{\sup}\notag \\
& &+\frac{d}{\delta}\sup_{s,t\in[S,1]}\left\|C(s,t)\right\|_{\sup}\left\|\tilde{\mathfrak{a}}(s_1,t_1)-\mathfrak{a}(s_2,t_2)\right\|_{\sup}\notag \\
& \leq& G^*\frac{d!d^2}{\delta^2}\sup_{s,t\in[S,1]}\|D(s,t)\|_{\sup}^{d-1}\sup_{s,t\in[S,1]}\left\|C(s,t)\right\|_{\sup}\sup_{s,t\in[S,1]}\left\|\tilde{\mathfrak{a}}(s,t)\right\|_{\sup}(\abs{s_1-s_2}+\abs{t_1-t_2})\notag \\
& &+G^*\frac{d!(d-1)}{\delta}\sup_{s,t\in[S,1]}\|D(s,t)\|_{\sup}^{d-2}\sup_{s,t\in[S,1]}\left\|\tilde{\mathfrak{a}}(s,t)\right\|_{\sup}(\abs{s_1-s_2}+\abs{t_1-t_2})\notag \\
& &+G^\prime\frac{d}{\delta}\sup_{s,t\in[S,1]}\left\|C(s,t)\right\|_{\sup}(\abs{s_1-s_2}+\abs{t_1+t_2})\notag \\
& =:&\tilde{G}[\abs{s_1-s_2}+\abs{t_1-t_2}],
}
and
\bqny{
\abs{\mathfrak{w}^\top(s_1,t_1)\tilde{\mathfrak{a}}(s_1,t_1)-\mathfrak{w}^\top(s_2,t_2)\tilde{\mathfrak{a}}(s_2,t_2)}&\leq& d\left\|\mathfrak{w}(s_1,t_1)-\mathfrak{w}(s_2,t_2)\right\|_{\sup}\left\|\tilde{\mathfrak{a}}(s_1,t_1)\right\|_{\sup}\\
&&+ d\left\|\mathfrak{w}(s_2,t_2)\right\|_{\sup}\left\|\tilde{\mathfrak{a}}(s_1,t_1)-\tilde{\mathfrak{a}}(s_2,t_2)\right\|_{\sup}\\
&\leq& d\tilde{G}\sup_{s,t\in[S,1]}\left\|\tilde{\mathfrak{a}}(s,t)\right\|_{\sup}(\abs{s_1-s_2}+\abs{t_1-t_2})\\
&&+ dG^\prime\sup_{s,t\in[S,1]}\left\|\mathfrak{w}(s,t)\right\|_{\sup}(\abs{s_1-s_2}+\abs{t_1+t_2})\\
&=:&\tilde{G}^\prime(\abs{s_1-s_2}+\abs{t_1-t_2}).
}
Since $\mathfrak{w}^\top(s,t)\tilde{\mathfrak{a}}(s,t)>0$ for all $s,t\in[S,1]$ we have that \eqref{holder} holds.
Applying Piterbarg inequality, we have that  there exist positive constants $C,\gamma$ such that
\bqny{
\pk{\exists s,t\in[S,1]:\mathcal{W}(s,t)-\mathfrak{c}(s,t)\geq\mathfrak{a}u} \leq Cu^{\gamma} e^{-(u+\mathfrak{C})^2/2\sigma^2}
}
for all $u$ positive.
Further, by Lemma \ref{proj1} for some constants $\alpha, C^*,C^+$ as $u\to\infty$
\bqny{
\lefteqn{
\pk{A_{\II\setminus\JJ}(s_1),A_{\JJ\setminus\II}(s_2),A_{\II\cap\JJ}(\min(s_1,s_2))}}\\
&=&\pk{\mathcal{W}(s_1,s_2)-\mathfrak{c}(s_1,s_2)\geq \mathfrak{a}u}\\
&\sim& C^*u^{-\alpha}e^{-\frac{1}{2}(\tilde{\mathfrak{a}}(s_1,s_2)u+
\mathfrak{c}(s_1,s_2))^\top \kk{D^{-1}(s_1,s_2)}(\tilde{\mathfrak{a}}(s_1,s_2)u+\mathfrak{c}(s_1,s_2))}\\
&=& C^+u^{-\alpha}e^{-\frac{u^2}{2\kk{\sigma^2}}}e^{-u(\tilde{\mathfrak{a}}_{s_1,s_2})^\top\kk{D^{-1}(s_1,s_2)}(\mathfrak{c}(s_1,s_2))}.
}
Hence the claim follows for $\tau=|\mathfrak{C}/\sigma^2|+\sup_{s,t\in[S,1]}|\tilde{\mathfrak{a}}(s,t)D^{-1}(s,t)\mathfrak{c}(s,t)|+1$.\QED

\begin{lem}\label{conv}
The function $r_\II(t),t>0$ defined in \eqref{rfunct} is convex and if $\vk c_\II$ has at least one positive component,
then there exists $T>0$ such that for some positive $s$ and any $t>0$
\bqn{ \label{claim1}
r_\II(T+t)\ge r_\II(T)+st.
}
Moreover, if $\vk a_\II+\vk c_\II t$ for any $t>0$ have at least one positive component, then $r_\II(t),t>0$ has a unique point of minimum.
\end{lem}

\prooflem{conv} For notational simplicity  we shall assume that $\II=\{1\ldot d\}$ and hence below omit in the notation
index $\II$.
Fix two points $0<t_1\leq t_2$, and some constant $\alpha\in(0,1)$. The convexity of $r$ follows by showing that
\bqny{
r(\alpha t_1+(1-\alpha)t_2)\leq\alpha r(t_1)+(1-\alpha)r(t_2).
}
Define below with  $\Sigma=\Gamma\Gamma^\top$, which is positive definite by the assumption on $\Gamma$
\bqny{
\vk x_1=arg\min_{\vk x\geq\vk a+\vk c t_1}\vk x^\top\Sigma^{-1}\vk x, \quad \vk x_2=arg\min_{\vk x\geq\vk a+\vk c t_2}\vk x^\top\Sigma^{-1}\vk x.
}
Since $\vk x_1\geq\vk a+\vk c t_1$ and $\vk x_2\geq\vk a+\vk c t_2$, then
\bqny{
\alpha \vk x_1+(1-\alpha)\vk x_2\geq\alpha(\vk a+\vk c t_1)+(1-\alpha)(\vk a+\vk c t_2)=\vk a+\vk c(\alpha t_1+(1-\alpha)t_2)
}
implying
$$r(\alpha t_1+(1-\alpha)t_2)\leq \frac{1}{\alpha t_1+(1-\alpha)t_2} (\alpha \vk x_1+(1-\alpha)\vk x_2)^\top\Sigma^{-1}(\alpha \vk x_1+(1-\alpha)\vk x_2).$$
Therefore it suffices to show that
\bqny{
\frac{1}{\alpha t_1+(1-\alpha)t_2} (\alpha \vk x_1+(1-\alpha)\vk x_2)^\top\Sigma^{-1}(\alpha \vk x_1+(1-\alpha)\vk x_2)\leq \frac{\alpha}{t_1} \vk x_1^\top\Sigma^{-1}\vk x_1
+\frac{1-\alpha}{t_2}\vk x_2^\top\Sigma^{-1}\vk x_2
}
i.e., we need to prove that
\bqny{
2\alpha(1-\alpha)\vk x_1\Sigma^{-1}\vk x_2\leq\alpha(1-\alpha)\left(\frac{t_2}{t_1}\vk x_1^\top\Sigma ^{-1}\vk x_1+\frac{t_1}{t_2}\vk x_2^\top\Sigma^{-1}\vk x_2\right).
}
Since $\alpha(1-\alpha)>0$, then
\bqny{
0\leq t_2^2\vk x_1^\top\Sigma^{-1}\vk x_1+t_1^2\vk x_2^\top\Sigma^{-1}\vk x_2-2t_1t_2\vk x_1^\top\Sigma^{-1}\vk x_2=(t_2\vk x_1-t_1\vk x_2)^\top\Sigma^{-1}(t_2\vk x_1-t_1\vk x_2),
}
which follows by the fact that $\Sigma^{-1}$ is positive definite.
Now, since $r$ is convex, $\lim_{t\to 0} r(t)\cE{=} \infty$
and $\Sigma^{-1}$ is positive definite,   there exists some positive vector
$\vk q\in\R^d$ such that $\Sigma^{-1}-Q$ is also positive definite,
where  $Q_{ij}=0$ if $i\not= j$ and $Q_{ii}=q_i$. Hence for some index $i$ such that $c_i>0$ and all $t$ such that $a_i+c_it>0$
\bqny{
r(t)\geq\frac{1}{t}\min_{x_i>a_{i}+c_i t}\vk x^\top\left(Q+\left(\Sigma^{-1}-Q\right)\right)\vk x\geq\frac{1}{t}q_i (a_{i}+c_i t)^2.
}
Consequently, $\limit{t} r(t)=\infty$ and hence $r(t)$ is monotone non-decreasing for all $t$ larger than some $T>0$,
implying the claim in \eqref{claim1}.

From the above arguments, $r(t),t>0$ is continuous, non-negative  with $\limit{t} r(t)= \IF$ and thus it attains its minimum on $(0,\IF)$.
Assume that its minimum is attained at two different points, say $t_1< t_2$. According to \nelem{prop1} for any $t>0$ there exists some index set $I_t$, such that
\bqny{
r(t)=\frac{1}{t}(\vk a_{I_t}+\vk c_{I_t} t)^\top(\Sigma_{I_tI_t})^{-1}(\vk a_{I_t}+\vk c_{I_t} t).
}
Clearly, there exists some uncountable set $Z\subset[t_1,t_2]$ such that $I_{t}$ is the same for all $t\in Z$. Moreover, as $r(t_1)=r(t_2)=\min_{t>0}r(t)$ and $r(t)$ is convex, $r(s)$ is the same for all $s\in Z$. But notice that
\bqny{
r(s)&=&\frac{1}{s}(\vk a_{I_s}+\vk c_{I_s} s)^\top(\Sigma_{I_sI_s})^{-1}(\vk a_{I_s}+\vk c_{I_s} s)\\
&=&\frac{1}{s}\vk a_{I_s}^\top(\Sigma_{I_sI_s})^{-1}\vk a_{I_s}+s\vk c_{I_s}^\top(\Sigma_{I_sI_s})^{-1}\vk c_{I_s}+2\vk a_{I_s}^\top(\Sigma_{I_sI_s})^{-1}\vk c_{I_s}.
}
Since this function is the same for any $s\in Z$ and $Z$ is uncountable, then it is constant, i.e., $\vk a_{I_s}=\vk 0_{I_s}$ and
$\vk{c}_{I_s}=\vk 0_{I_s}$. In this case $r(t)=0$. Since $\Sigma$ is positive definite, then $\vk a+\vk c s=0$ for all $s\in Z$, which contradicts with our assumption. Hence the claim follows.
\QED

\begin{lem}{\label{nonstat}}
Suppose that  $\Sigma=\Gamma\Gamma^\top$ is positive definite. For any non-empty subset $\II\subset\{1\ldot d\}$ if $\vk c_\II$ and $\vk a_\II+\vk c_\II t$ for all $t\geq 0$ have at least one positive component, then for any point $0<t\not=\hat{t}=arg\min_{t>0} r_\II(t)$ there exists some positive constant $\nu$ such that
\bqny{
\pk{\vk W_\II(t)>(\vk a_\II+\vk c_\II t)\sqrt{u}}=o\left(e^{-\nu u}\right)\pk{\vk W_\II(\hat{t})>(\vk a_\II+\vk c_\II\hat{t})\sqrt{u}}, \quad u\to \IF.
}
\end{lem}
\prooflem{nonstat} For notational simplicity we omit below the subscript $\II$. Since for any $t>0$ we have $Var(\vk W(t))=t\Sigma$, then by Lemma \ref{proj1}
\bqny{
\pk{\vk W(t)>(\vk a+\vk ct)\sqrt{u}}\sim Cu^{-\alpha(t)/2}e^{-\frac{u}{2t}\tilde{\vk{\mathfrak{p}}}(t)^\top\Sigma^{-1}\tilde{\vk{\mathfrak{p}}}(t)},
}
where $C$ is some positive constant, $\alpha(t)$ is an integer and $\tilde{\vk{\mathfrak{p}}}(t)$ is the unique solution of $\Pi_{t\Sigma}(\vk a+\vk ct)$, which can be reformulated also as
\bqny{
\pk{\vk W(t)>(\vk a+\vk ct)\sqrt{u}}\sim Cu^{-\alpha(t)/2}e^{-\frac{u}{2}r(t)}, \quad u\to \IF.
}
If $t\not=\hat{t}$, then $r(t)-r(\hat{t})=\tau>0$ and
\bqny{
\frac{\pk{\vk W(t)>(\vk a+\vk ct)\sqrt{u}}}{\pk{\vk W(\hat{t})>(\vk a+\vk c\hat{t})\sqrt{u}}}\sim C^*u^{(\alpha(\hat{t})-\alpha(t))/2}e^{-\frac{\tau u}{2}}=o\left(e^{-\frac{\tau}{3}u}\right)
}
as $u\to \IF$.
\QED

\begin{lem}{\label{pit}} Let $\vk a,\vk c\in\R^d$ be such that $\vk a+\vk c t$ has at least one positive component for all $t$
in a compact set $\mathcal{T}\subset (0, \infty)$. If  $\Sigma=\Gamma\Gamma^\top$ is positive definite, then there exist constants $C>0$, $\gamma>0$  and $\mathfrak{t}\in\mathcal{T}$ such that
for all $u>0$
\bqny{
\pk{\exists t\in\mathcal{T}:~\vk W(t)>(\vk a+\vk c t)\sqrt{u}}\leq Cu^{\gamma}e^{-\frac{u}{2}r(\mathfrak{t})}.
}
If we also have that for some non-overlapping index sets $\II,\JJ\subset\{1\ldot d\}$ and some compact subset
$\mathcal{T}\subset [0,\infty)^2$   both
$((\vk a_{\II}+\vk c_{\II} t_1)^\top, (\vk a_{\JJ}+\vk c_{\JJ} t_2)^\top)^\top$
\kk{have} at least one positive component for all $(t_1,t_2)\in\mathcal{T}$, then for some $\mathfrak{t}=(\mathfrak{t}_1,\mathfrak{t}_2)\in\mathcal{T}$   as $u\to\infty$
\bqny{
& \pk{\exists t\in\mathcal{T}:~\vk W_\II(t_1)>(\vk a_\II+\vk c_\II t_1)\sqrt{u},~\vk W_\JJ(t_{\kk{2}})>(\vk a_\JJ+\vk c_\JJ t_2)\sqrt{u}}\\
& \qquad\qquad= o(e^{\sqrt{u}}\pk{\vk W_\II(\kk{\mathfrak{t}}_1)>(\vk a_\II+\vk c_\II \mathfrak{t}_1)\sqrt{u},
~\vk W_\JJ(\kk{\mathfrak{t}_2})>(\vk a_\JJ+\vk c_\JJ \kk{\mathfrak{t}_2})\sqrt{u}}).
}
\Ee{Moreover, the same estimate holds} if $\II$ and $\JJ$ are overlapping and for all $(\mathfrak{t}_1,\mathfrak{t_2})\in\mathcal{T}$ we have $\mathfrak{t}_1\not=\mathfrak{t}_2$.
\end{lem}

\prooflem{pit} Denote by $D(t)$ the covariance matrix of
$\vk W(t)$, which by assumption on $\Gamma$ is positive definite.
Let $\tilde{\mathfrak{a}}(t) = arg\min\limits_{\vk x\geq\vk a+\vk c t} \vk x^\top D^{-1}(t) \vk x$
be the solution of $\Pi_{D}(\vk a+\vk c t), t>0$ and let further
\bqny{
\mathfrak{w}(t) =D ^{-1}(t)\tilde{\mathfrak{a}}(t)
}
be the solution of the dual optimization problem. In view of \eqref{MC}  $\mathfrak{w}_I(t)$ has positive  components \COM{for $I$ the unique index set related to
$\Pi_{D(t)}(\vk a+\vk c t)$ } and moreover by \eqref{eq:new}
$$f(t)=\mathfrak{w}^\top(t) (\vk a+ \vk c t) = \tilde{\mathfrak{a}}^\top(t) D^{-1}(t)\tilde{\mathfrak{a}}(t)>0
$$
implying
\bqny{
\pk{\exists t\in\mathcal{T}:\vk W(t)\geq(\vk a+\vk c t) \sqrt{u}}&\leq&
\pk{\exists t\in\mathcal{T}:\frac{\mathfrak{w}^\top \kk{(t)} \vk W(t)}{\mathfrak{w}^\top \kk{(t)}(\vk a+\vk c t)}\geq \sqrt{u}}.
}
We have further that
\bqny{
\sigma^{2}=\sup_{t\in\mathcal{T}}\E{\left(\frac{\mathfrak{w}^\top(t)\vk W(t)}{\mathfrak{w}^\top(t) (\vk a+\vk c t)}\right)^2}
=\sup_{t\in\mathcal{T}}\frac{1}{\tilde{\mathfrak{a}}^\top(t) D^{-1}(t)\tilde{\mathfrak{a}}(t)}
=\frac{1}{\tilde{\mathfrak{a}}^\top(\mathfrak{t})D^{-1}(\mathfrak{t})\tilde{\mathfrak{a}}(\mathfrak{t})}>0
}
for some $\mathfrak{t}\in\mathcal{T}$, since $\mathcal{T}$ is compact. Since $f(t)>0, t\in\mathcal{T}$ is continuous,  we may apply Piterbarg inequality (\Ee{as in the proof of \eqref{LL4.1}}) and obtain
\bqny{
\pk{\exists t\in\mathcal{T}:\vk W(t)\geq(\vk a+\vk ct)\sqrt{u}} \leq Cu^{\gamma} e^{-u/2\sigma^2}
}
for some positive constants $\gamma$ and $C$, which depend only on $\vk W(t)$ and $d$.
Since, by the definition we have $r(\mathfrak{t})=1/\sigma^2$,
the proof of the first inequality is complete.

The next assertion may be obtained with the same arguments but for vector-valued random process
\bqny{
\mathcal{W}(s,t)=(\vk W_{\II}^\top(s),\vk W_{\JJ}^\top(t))^\top.
}
By the definition of $\mathcal{T}$, for any $(s,t)\in\mathcal{T}$ we have
$\lvert Var(\mathcal{W}(s,t))\rvert >0$, thus we can apply Piterbarg inequality and
in consequence, using \nelem{proj1}, the claim follows.
\QED

\begin{lem}{\label{rtail2}} Suppose that  $\Sigma=\Gamma\Gamma^\top$ is positive definite. For any subset $\II\subset \{1\ldot d\}$ if $\vk c_\II \in \R^{\abs{\II}}$ has at least one positive component and $\vk a_\II+\vk c_\II t \inr^{\abs{\II}}$ has at least one positive component for all non-negative $t$, then for some positive constants $\nu$, $\hat{t}=arg\min\limits_{t>0} r_\II(t)$ and all $T$ large
\bqny{
\pk{\exists t>T:~\vk W_\II(t)>(\vk a_\II+\vk c_\II t)\sqrt{u}}=o(e^{-\nu \kk{u}})
\pk{\vk W_\II(\hat{t})>(\vk a_\II+\vk c_\II \hat{t})\kk{\sqrt{u}}},\qquad u\to\infty.
}
\end{lem}
\prooflem{rtail2} For notational simplicity we omit below the subscript  $\II$.
For some given $T>\hat{t} $ we have using Lemmas \ref{pit}, \ref{conv}
\bqny{
\pk{\exists t>T:~\vk W(t)>(\vk a+\vk ct)\sqrt{u}} &\leq& \sum_{i=0}^{\infty}\pk{\exists t\in [T+i,T+i+1]:~\vk W(t)>(\vk a+\vk ct)\sqrt{u}}\\
&\leq&\sum_{i=0}^{\infty}Cu^{\gamma}e^{-\frac{r(\mathfrak{t}_i)}{2}u}\\
&\leq & Cu^{\gamma}e^{-\frac{r(T)}{2}u}\sum_{i=0}^{\infty}e^{-isu}
\\
&\leq& Cu^{\gamma}e^{-\frac{r(T)}{2}u}\left(1+\int_0^\infty e^{-sux}\td x\right),
}
where the last integral is finite and decreasing for sufficient large $u$. Taking thus $u>1$ we obtain
\bqny{
\pk{\exists t>T:~\vk W(t)>(\vk a+\vk ct)u}\leq C^*u^{\gamma}e^{-\frac{r(T)}{2}u^2}
}
for some $C^*>0$. Hence the claim follows with the same arguments as in the proof of Lemma \ref{nonstat}.
\QED

\prooflem{diff} Using Lemma \ref{rtail2} we know that there exist points $t_\II,~t_\JJ$ such that
\bqny{
\pk{\exists t\geq T_\II:~A^*_\II(t)}=
o(\pk{A^*_\II(\hat{t}_\II)}),\qquad \pk{\exists t\geq T_\JJ:~A^*_\JJ(t)}=o(\pk{A^*_\JJ(\hat{t}_\JJ)}),\qquad u\to\infty.
}
Next, for some   positive $\varepsilon<|\hat{t}_\II-\hat{t}_\JJ|/3$ we have
\bqny{
\pk{\exists s,t>0:~A^*_{\II}(t) \kk{\cap} A^*_{\JJ}(s)}&\leq&\pk{\exists (s,t)\in[\hat{t}_\II-\varepsilon,\hat{t}_\II+\varepsilon]\times[\hat{t}_\JJ-\varepsilon,\hat{t}_\JJ+\varepsilon]:~A^*_{\II}(t)\cap A^*_{\JJ}(s)}\\
&+&\pk{\exists t\in[0,\hat{t}_\II-\varepsilon]:~A^*_{\II}(t)}+\pk{\exists t\in[\hat{t}_\II+\varepsilon,T_\II]:~A^*_{\II}(t)}\\
&+&\pk{\exists t\in[0,\hat{t}_\JJ-\varepsilon]:~A^*_{\JJ}(t)}+\pk{\exists t\in[\hat{t}_\JJ+\varepsilon,T_\JJ]:~A^*_{\JJ}(t)}\\
&+&\pk{\exists t\geq T_{\II}:~A^*_{\II}(t)}+\pk{\exists t\geq T_{\JJ}:~A^*_{\JJ}(t)}.
}
Using Lemmas \ref{pit}, \ref{rtail2} and
$$\pk{A^*_\II(t)}\sim Cu^{-\alpha}e^{-r(t)u/2}, \quad \pk{A^*_\II(t)}=o(ue^{-r(t)u/2}), \quad u\to \IF$$
we obtain
\bqny{
\pk{\exists s,t>0:~A^*_{\II}(t)\kk{\cap} A^*_{\JJ}(s)}
&\kk{=}& o(e^{\sqrt{u}} \pk{A^*_{\II}(s_1)\cap A^*_{\JJ}(s_2)})\\
&+&o(u^{\tau_3}\pk{A^*_{\II}(t_3)})+o(u^{\tau_4}\pk{A^*_{\II}(t_4)})\\
&+&o(u^{\tau_5}\pk{A^*_{\JJ}(t_5)})+o(u^{\tau_6}\pk{A^*_{\JJ}(t_6)})\\
&+&o(\pk{A^*_\II(\hat{t}_{\II})})+o(\pk{A^*_{\JJ}(\hat{t}_{\JJ})})
}
for some positive constants \kdd{$t_i, 3 \le i \le 6$,} where
\bqny{
t_3\in[0,\hat{t}_\II-\varepsilon],\quad t_4\in[\hat{t}_\II+\varepsilon,T_\II], \quad t_5\in[0,\hat{t}_\JJ-\varepsilon],\quad t_6\in[\hat{t}_\JJ+\varepsilon,T_\JJ]\quad s_1\in[\hat{t}_\II-\varepsilon,\hat{t}_\II+\varepsilon]\quad s_2\in[\hat{t}_\JJ-\varepsilon,\hat{t}_\JJ+\varepsilon].
}
Note that for $i=3,4$, $t_i\not=\hat{t}_{\II}$. Hence by Lemma \ref{nonstat}
\bqny{
u^{\tau_i}\pk{A^*_{\II}(t_i)}=o(\pk{A^*_\II(\hat{t}_{\II})}).
}
The same works also for $j=5,6$
\bqny{
u^{\tau_j}\pk{A^*_{\JJ}(t_j)}=o(\pk{A^*_\JJ(\hat{t}_{\JJ})}).
}
Thus we can focus only on the first probability.
By the definition of $A^*_\II$ and $A^*_\JJ$ in \Ee{\eqref{AY}}
\bqny{
\pk{A^*_\II(s_1)\cap A^*_\JJ(s_2)}=\pk{ \mathcal{W}(s_1,s_2) > \vk b \sqrt{u}},
}
where $\vk b= ((\vk a_\II+\vk c_\II s_1)^\top,(\vk a_\JJ+\vk c_\JJ s_2)^\top)$ and
\bqny{
\mathcal{W}(s,t)=(\vk W_{\II}(s)^{\top},\vk W_{\JJ}(t)^\top)^\top.
}
 Define $\widehat{i}=\II\cup\JJ\setminus \{i\}$. Applying Remark \ref{reduct}, there exists an index $i$ and a constant $\eta>0$ such that
\bqny{
\pk{A^*_\II(s_1)\cap A^*_\JJ(s_2)}=o\left(e^{-\eta u}\right)\pk{ (\mathcal{W}(s_1,s_2))_{\widehat{i}} > \vk b_{\widehat{i}} \sqrt{u}}.
}
If $i\in\II$, then
\bqny{
\pk{ (\mathcal{W}(s_1,s_2))_{\widehat{i}}>\vk b_{\widehat{i}}\sqrt{u}}\leq\pk{\vk W_\JJ(s_2)>(\vk a_\JJ+\vk c_\JJ s_2)u},
}
\Ee{or}
\bqny{
\pk{(\mathcal{W}(s_1,s_2))_{\widehat{i}} > \vk b_{\widehat{i}} \sqrt{u}}\leq\pk{\vk W_\II(s_1)>(\vk a_\II+\vk c_\II s_1)u}.
}
In both cases
\bqny{
e^{\sqrt{u}}\pk{A^*_\II(s_1)\cap A^*_\JJ(s_2)}&=&o\left(\pk{\vk W_\II(s_1)>(\vk a_\II+\vk c_\II s_1)u}+\pk{\vk W_\JJ(s_1)>(\vk a_\JJ+\vk c_\JJ s_1)u}\right)\\
&=&o\left(\pk{\vk W_\II(\hat{t}_\II)>(\vk a_\II+\vk c_\II \hat{t}_\II)u}+\pk{\vk W_\JJ(\hat{t}_\II)>(\vk a_\JJ+\vk c_\JJ \hat{t}_\JJ)u}\right)
}
establishing the proof.
\QED

\prooflem{same} Using Lemma \ref{rtail2} we have
\bqny{
\pk{\exists s,t>0:~A^*_\II(s)\kk{\cap} A^*_\JJ(t)}&\leq& \pk{\exists (s,t)\in\mathbb{T}_1:~A^*_\II(s)\kk{\cap} A^*_\JJ(t)}+\pk{\exists (s,t)\in\mathbb{T}_2:~A^*_\II(s)\kk{\cap} A^*_\JJ(t)}\\
&+&o(\pk{A^*_\II(\hat{t}_\II)})+o(\pk{A^*_\JJ(\hat{t}_\JJ)}),
}
where
\bqny{
\mathbb{T}_1=\{(s,t)\in[0,T_\II]\times[0,T_\JJ]:|s-\hat{t}_\II|\geq|t-\hat{t}_\II|\},\\
\mathbb{T}_2=\{(s,t)\in[0,T_\II]\times[0,T_\JJ]:|s-\hat{t}_\II|\leq|t-\hat{t}_\II|\}
}
and $T_{\II}$ and $T_\JJ$ are the constants from \eqref{claim1}.
According to \nelem{pit} for some $(s_i,t_i)\in\mathbb{T}_i$
\bqny{
\pk{\exists (s,t)\in\mathbb{T}_i:A^*_\II(s)\cap A^*_\JJ(t)}=o \Bigl(e^{\sqrt{u}} \Bigr) \pk{A^*_\II(s_i)\cap A^*_{\JJ\setminus\II}(t_i)}.
}
If $s_1\not=\hat{t}_{\II}$, then according to Lemma \ref{nonstat}
\bqny{
e^{\sqrt{u}}\pk{A^*_\II(s_1)\cap A^*_{\JJ\setminus\II}(t_1)}\leq e^{\sqrt{u}}\pk{A^*_\II(s_1)}=o(\pk{A^*_{\II}(\hat{t}_{\II})}).
}
Otherwise, using the definition of \kk{$\mathbb{T}_1$}, $|t_1-\hat{t}_\II|\leq|s_1-\hat{t}_{\II}|=0$, so $t_1=\hat{t}_{\II}$ and thus
\bqny{
\pk{A^*_\II(s_1)\cap A^*_{\JJ\setminus\II}(t_1)}=\pk{A^*_{\II\cup\JJ}(\hat{t}_{\II})}.
}
This probability can be bounded using Remark \ref{reduct}, namely we have
\bqny{
\pk{A^*_{\II\cup\JJ}(\hat{t}_{\II})}=o\Bigl(e^{-\nu \kk{u}} \Bigr)\pk{A^*_{\II\cup\JJ\setminus\{i\}}(\hat{t}_{\II})}
}
for some $i\in\II\cup\JJ$ and $\eta>0$. As $|\II|=|\JJ|=k$, and $\II\not=\JJ$, then $|\II\cup\JJ|\ge k+1$ and thus $|\II\cup\JJ\setminus\{i\}|\ge k$.
Consequently, we have
\bqny{
e^{\sqrt{u}}\pk{A^*_{\II\cup\JJ}(\hat{t}_{\II})}=o\Bigl(\pk{A^*_{\II\cup\JJ\setminus\{i\}}(\hat{t}_{\II})} )=o\left(\sum_{\substack{\KK\subset\{1\ldots d\}\\|\KK|=k}}\pk{A^*_\KK(\hat{t}_{\KK})}\right).
}

With similar arguments we obtain further
\bqny{
\pk{\exists (s,t)\in\mathbb{T}_2:~A^*_\II(s)\kk{\cap} A^*_\JJ(t)}=o\left(\sum_{\substack{\KK\subset\{1\ldots d\}\\|\KK|=k}}\pk{A^*_\KK(\hat{t}_{\KK})}\right).
}
Hence the claim follows.

\QED

\Ee{Recall that $\tilde{ \vk a}$ stands for} the unique solution of the quadratic programming problem $\Pi_\Sigma(\vk a)$.

\prooflem{inimp:Lem}
By the self-similarity of Brownian motion \Ee{for all $u>0$}
\bqny{
m(u,\CN{\Lambda}):=\pk{\exists_{t\in[0,\delta(u,\CN{\Lambda})]}:\vk W(t)-t\vk c> u\vk a}=\pk{\exists_{t\in[0,1]}:\vk W(t)-\delta^{1/2}(u,\CN{\Lambda})t\vk c> \delta^{-1/2}(u,\CN{\Lambda})u\vk a}.
}
Applying  \kk{Theorem} \ref{D} we obtain
\bqny{
m(u,\CN{\Lambda})\le\frac{\pk{\vk W(1)\ge \delta^{-1/2}(u,\CN{\Lambda})u\vk a +\delta^{1/2}(u,\CN{\Lambda})\vk c}}{\pk{\vk W(1)>\max(\vk c,\vk 0)}}.
}
Further, for all $u$ large
\bqny{
\lefteqn{\pk{\vk W(1)\ge \delta^{-1/2}(u,\CN{\Lambda})u\vk a +\delta^{1/2}(u,\CN{\Lambda})\vk c}}\\
&\leq&
\pk{\vk W(1)\ge \delta^{-1/2}(u,\CN{\Lambda})(u\vk a + \delta(u,\CN{\Lambda})\vk c)}\\
&=&
\int_{ \vk x \ge \delta^{-1/2}(u,\CN{\Lambda})(u\vk a+\delta(u,\CN{\Lambda})\vk c)}\varphi(\vk x)\td\vk x\\
&=&\frac{1}{\delta^{d/2}(u,\CN{\Lambda})}\int_{ \vk x \ge u\vk a+\delta(u,\CN{\Lambda})\vk c }\varphi\left(\delta^{-1/2}(u,\CN{\Lambda})\vk x\right)\td\vk x\\
&=&\frac{1}{\delta^{d/2}(u,\CN{\Lambda})}\int_{ \vk x \ge u\vk a+\vk c }\varphi\left(\delta^{-1/2}(u,\CN{\Lambda}/2)(\vk x-(1-\delta(u,\CN{\Lambda}))\vk c)\right)\td\vk x\\
&=&\frac{1}{\delta^{d/2}(u,\CN{\Lambda})}\int_{ \vk x \ge u\vk a+\vk c }\varphi\left(\delta^{-1/2}(u,\CN{\Lambda}/2)(\vk x-\CN{\Lambda}u^{-2}\vk c)\right)\td\vk x,
}
where $\varphi$ is the pdf of $\vk W(1)$.
Next, we bound the integrand above considering $\vk x\ge  u \vk a+ \vk c$ as follows
\bqny{
\lefteqn{\varphi\left(\delta^{-1/2}(u,\CN{\Lambda})(\vk x-\CN{\Lambda}u^{-2}\vk c)\right)}\\
&=&\frac{1}{(2\pi)^{d/2}\kk{\sqrt{|\Sigma|}}}\exp\left(-\frac{1}{2\delta(u,\CN{\Lambda})}(\vk x-\CN{\Lambda}u^{-2}\vk c)^\CN{\top}\Sigma^{-1}(\vk x-\CN{\Lambda}u^{-2}\vk c)\right)\\
&\le &
\frac{1}{(2\pi)^{d/2}\kk{\sqrt{|\Sigma|}}}\exp\left(-\frac{1+\CN{\Lambda}u^{-2}}{2}(\vk x-\CN{\Lambda}u^{-2}\vk c)^\CN{\top}\Sigma^{-1}(\vk x-\CN{\Lambda}u^{-2}\vk c)\right)\\
&=&
\frac{1}{(2\pi)^{d/2}\sqrt{|\Sigma|}}\exp\left(-\frac{1}{2}(\vk x-\CN{\Lambda}u^{-2}\vk c)^\CN{\top}\Sigma^{-1}(\vk x-\CN{\Lambda}u^{-2}\vk c)\right)\\
&\times&\exp\left(-\frac{\CN{\Lambda}u^{-2}}{2}(\vk x-\CN{\Lambda}u^{-2}\vk c)^\CN{\top}\Sigma^{-1}(\vk x-\CN{\Lambda}u^{-2}\vk c)\right)\\
&=&
\frac{1}{(2\pi)^{d/2}\sqrt{|\Sigma|}}\exp\left(-\frac{1}{2}(\vk x-\CN{\Lambda}u^{-2}\vk c)^\CN{\top}\Sigma^{-1}(\vk x-\CN{\Lambda}u^{-2}\vk c)\right)\\
&\times&\exp\left(-\frac{\CN{\Lambda}u^{-2}}{2}\left(\frac{\vk x}{\sqrt{2}}-\sqrt{2}\CN{\Lambda}u^{-2}\vk c\right)^\CN{\top}\Sigma^{-1}\left(\frac{\vk x}{\sqrt{2}}-\sqrt{2}\CN{\Lambda}u^{-2}\vk c\right)\right)\\
&\times&\exp\left(-\frac{\CN{\Lambda}u^{-2}}{4}\vk x^\top\Sigma^{-1}\vk x\right)\exp\left(\frac{\CN{\Lambda}^3u^{-6}}{4}\vk c^\top\Sigma^{-1}\vk c\right)\\
&\leq&
\frac{1}{(2\pi)^{d/2}\sqrt{|\Sigma|}}\exp\left(-\frac{1}{2}(\vk x-\CN{\Lambda}u^{-2}\vk c)^\CN{\top}\Sigma^{-1}(\vk x-\CN{\Lambda}u^{-2}\vk c)\right)\\
&\times&\exp\left(-\frac{\CN{\Lambda}}{4}(\vk x/u)^\top\Sigma^{-1}(\vk x/u)\right)\exp\left(\frac{\CN{\Lambda}^3u^{-6}}{4}\vk c^\top\Sigma^{-1}\vk c\right)\\
&\leq&
\varphi(\vk x)\exp\left(\frac{1}{2}(\sqrt{\CN{\Lambda}}u^{-3/2}\vk x)^\top\Sigma^{-1}(\sqrt{\CN{\Lambda}}u^{-3/2}\vk x)+\frac{1}{2}(u^{-1/2}\vk c)^\top\Sigma^{-1}(u^{-1/2}\vk c)\right)\\
&\times&\exp\left(-\frac{\CN{\Lambda}^2u^{-4}}{2}\vk c^\top\Sigma^{-1}\vk c\right)\exp\left(-\frac{\CN{\Lambda}}{4}(\vk x/u)^\top\Sigma^{-1}(\vk x/u)\right)\exp\left(\frac{\CN{\Lambda}^3u^{-6}}{4}\vk c^\top\Sigma^{-1}\vk c\right)\\
&\leq&
\varphi(\vk x)\exp\left(\left(-\frac{\CN{\Lambda}^2u^{-4}+u^{-1/2}}{2}+\frac{\CN{\Lambda}^3u^{-6}}{4}\right)\vk c^\top\Sigma^{-1}\vk c\right)\\
&\times&\exp\left(-\frac{\CN{\Lambda}-2\CN{\Lambda}u^{-1/2}}{4}(\vk x/u)^\top\Sigma^{-1}(\vk x/u)\right).
}

For large enough $u$, we have that $\vk x/u\geq \vk a+\vk c/u>\vk a+\varepsilon\min(\vk c,\vk 0)$ for all $\varepsilon>0$
\kdd{and}
\bqny{
C^*:=\min_{\vk y\geq \vk a+\varepsilon\min(\vk c,0)}\vk y^\top\Sigma^{-1}\vk y>0.
}
Hence, for large enough $u$
\bqny{
\exp\left(-\frac{\CN{\Lambda}-2\CN{\Lambda}u^{-1/2}}{4}(\vk x/u)^\top\Sigma^{-1}(\vk x/u)\right)\leq\exp\left(\frac{-C^*\CN{\Lambda}}{8}\right),\\
\exp\left(\left(-\frac{\CN{\Lambda}^2u^{-4}+u^{-1/2}}{2}+\frac{\CN{\Lambda}^3u^{-6}}{4}\right)\vk c^\top\Sigma^{-1}\vk c\right)\leq \exp\left(\frac{C^*\CN{\Lambda}}{16}\right),
}
\kdd{which imply that, for sufficiently large $u$}
\bqny{
\varphi\left(\delta^{-1/2}(u,\CN{\Lambda})(\vk x-\CN{\Lambda}u^{-2}\vk c)\right)\leq\varphi(\vk x)\exp\left(\frac{-\CN{\Lambda}}{16/C^*}\right), \quad \vk x\ge u \vk a + \vk c.
}

Consequently,  for sufficient large $u$ and \Ee{some} positive $C$
\bqny{
m(u,\CN{\Lambda}) &\leq& \frac{1}{\delta^{d/2}(u,\CN{\Lambda}/2)}\frac{\int_{ \vk x \ge u\vk a+\vk c }\varphi\left(\delta^{-1/2}(u,\CN{\Lambda}/2)(\vk x-\CN{\Lambda}u^{-2}\vk c)\right)\td\vk x}{\pk{\vk W(1)>\max(\vk c, \vk 0)}}\\
&\leq& \frac{\exp\left(\frac{-\CN{\Lambda}}{16/C^*}\right)}{\delta^{d/2}(u,\CN{\Lambda}/2)}\frac{\int_{\vk x \ge u\vk a+\vk c}\varphi(\vk x) \td \vk x}{\pk{\vk W(1)>\max(\vk c,\vk 0)}}\\
&=& \frac{\exp\left(\frac{-\CN{\Lambda}}{16/C^*}\right)}{\delta^{d/2}(u,\CN{\Lambda}/2)}\frac{\pk{\vk W(1)>\vk au+\vk c}}{\pk{\vk W(1)>\max(\vk c,\vk 0)}}\\
&\leq& \exp\left(-\frac{\CN{\Lambda}}{C}\right)\frac{\pk{\vk W(1)>\vk au+\vk c}}{\pk{\vk W(1)>\max(\vk c,\vk 0)}},
}
where $\vk 0\in\R^d$ has all components equal to 0. Hence \CN{\eqref{inimp}} follows.

\Ee{We show next \eqref{A}.}
Let   $I$ and $J$ are be the unique index sets corresponding to the quadratic programming problems $\Pi_{\Sigma}(\vk a).$
Define the vector
\begin{eqnarray*}
\vk v=u\tilde{\vk a}+\vk c- \vk x / \bar{\vk u},
\end{eqnarray*}
where $\bar{\vk u}$ has the components in the set $I$ equal to $u$ and the other components equal to 1 and
$\tilde{\vk a}$ is the unique solution of $\Pi_{\Sigma}(\vk a).$
Recall that we have the representation
\begin{eqnarray*}
\vk W(t)= \Gamma \vk B(t), \quad t\ge 0,
\end{eqnarray*}
where $\vk B(t)=(B_1(t)\ldot B_d(t))^\top$ with $B_i$'s independent standard Brownian motions
and \kk{${\Gamma}$} a $d\times $d real-valued non-singular matrix. Writing $\varphi$ for the pdf of $\vk W(1)$ we have
\begin{eqnarray*}
\varphi(\vk v)=\frac{1}{(2\pi)^{d/2}\sqrt{|\Sigma|}}e^{-\frac{1}{2}\left(u\tilde{\vk a}+\vk c-\vk x/\bar{\vk u}\right)^\top\Sigma^{-1}\left(u\tilde{\vk a}+\vk c-\vk x/\bar{\vk u}\right)}=\varphi(u\tilde{\vk a} +\vk c)\theta_u(\vk x)
\end{eqnarray*}
for all $\vk x\inr^d, u>0$,
where
\begin{eqnarray*}
\theta_u(\vk x)=e^{u\tilde{\vk a}^\top\Sigma^{-1}(\vk x/\bar{\vk u})}e^{\vk c^\top\Sigma^{-1}(\vk x/\bar{\vk u})}e^{-\frac{1}{2}(\vk x/\bar{\vk u})^\top \Sigma^{-1}(\vk x/\bar{\vk u})}.
\end{eqnarray*}
Letting
$\vk \lambda= \Sigma^{-1} \tilde{\vk a}$ by \eqref{MC}
$\vk \lambda_I= (\Sigma_{II})^{-1}\vk a_I> \vk 0_I $
and further \Ee{as $u\to \IF$}
\bqny{
& &e^{-\frac{1}{2}(\vk x/\bar{\vk u})^\top \Sigma^{-1}(\vk x/\bar{\vk u})}\to e^{-\frac{1}{2}\vk x_J^\top(\Sigma^{-1})_{JJ}\vk x_J},\\
& &e^{\vk c^\top\Sigma^{-1}(\vk x/\bar{\vk u})}\to e^{(\vk c_I^\top(\Sigma^{-1})_{IJ}+\vk c_J^\top(\Sigma^{-1})_{JJ})\vk x_J}, \ e^{\tilde{\vk a}^\top\Sigma^{-1}(u\vk x/\bar{\vk u})}\to e^{\vk a_I^\top(\Sigma^{-1})_{II}\x_I}= e^{\vk \lambda_I^\top \vk x_I}.
}
For $\tilde{\vk c}$  such that
$$(\vk c_I^\top(\Sigma^{-1})_{IJ}+\vk c_J^\top(\Sigma^{-1})_{JJ})=\tilde{\vk c}(\Sigma^{-1})_{JJ}$$
we have
\bqn{
\limit{u} \theta_{u}(\vk x) = \theta(\vk x)=e^{-\frac{1}{2}(\vk x_J-\tilde{\vk c})^\top(\Sigma^{-1})_{JJ}(\vk x_J-\tilde{\vk c})}e^{\frac{1}{2}\tilde{\vk c}^\top(\Sigma^{-1})_{JJ}\tilde{\vk c}}e^{\vk \lambda_I^\top \vk x_I}.
\label{psiasympt}
}
Setting
$F(u)=u^{-|I|}\varphi(u\tilde{\vk a}+\vk c)$
and
\bqny{
h_u([0,\CN{\Lambda}],\vk x)=\pk{\exists_{t\in[0,\CN{\Lambda}]}:\Gamma\vk B(\overline{t})-\overline{t}\vk c>u\vk a|\vk B(1)=  \Gamma ^{-1}\vk v}
}
we may write further
\bqny{
M(u,\CN{\Lambda})&=&u^{-|I|}\int_{\R^d}\pk{\exists_{t\in[\delta(u,\CN{\Lambda}),1]}:\vk W(t)-t\vk c>u\vk a| \vk W(1)=\vk v}\varphi(\vk v)\td\vk x\\
&=&F(u)\int_{\R^d} h_u([0,\CN{\Lambda}],\vk x)\theta_u(\vk x)\td\vk x,
}
where $\overline{t}=1-t/u^2$. Define
\bqny{
\vk B_{u}(t)+\Gamma ^{-1}\vk v \overline{t}:= \vk B(\overline{t})|\vk B(1)=\Gamma ^{-1}\vk v
}
and note that  for  all $t,u$ positive
\bqny{
  u( \Gamma \vk B(\overline{t})-\overline{t}\vk c-u\vk a)\Bigl \lvert \vk W(1)=\vk v &=& u(\Gamma \vk B_u(t)+\overline{t}\vk v -\overline{t}\vk c-u\vk a)\\
&=& u(\Gamma \vk B_u(t)+\left(1-t/u^2\right)(u\tilde{\vk a}-\vk x/\bar{\vk u})-u\vk a)\\
&=&\Gamma (u\vk B_{u}(t))+u^2(\tilde{\vk a}-\vk a)-u(\vk 0_I^\top,\vk x_J^{\top})^{\top}-(\vk x_I^\top,\vk 0_J^\top)^\top-t\tilde{\vk a}+u^{-1}t\vk x/\bar{\vk u}.
}
According to this representation, with similar arguments as in \cite{mi:18}, using continuous mapping theorem
we obtain for almost all $\vk x$ (with respect to the Lebesgues measure on $\R^d$)
\bqny{
\limit{u} h_{u}([0,\CN{\Lambda}],\vk x) = h([0,\CN{\Lambda}],\vk x),
}
where
\bqn{
h([0,\CN{\Lambda}],\vk x)&=&
\pk{\exists_{t\in[0,\CN{\Lambda}]}:\vk W_{I }(t)-\vk x_I-t\tilde{\vk a}_{I }>\vk 0_{I}}\mathbb{I}_{\{\vk x_{U}<\vk 0_{U}\}},\\
U&=&\{j\in J:\tilde{ a}_j= a_j\}.
}
For \EZ{notational simplicity}  we assume next that $U$ defined above  is non-empty.
In order to use the dominated convergence theorem, we need to find bounds for the functions $\theta_u(\vk x)$ and $h_u([0,T],\vk x)$.\\
Fix some $\varepsilon\in\left(0,\min_{i\in I}\lambda_i\right)$ such that the matrix $\Sigma^{-1}-\varepsilon I_d$ is positive definite (here $I_d$ is the identity matrix). We write next $< \cdot , \cdot>$ for the scalar product on $\R^l, l\le d$. For any $u>0$ we have
\bqny{
\theta_u(\vk x)&=&e^{<\vk\lambda_I,\vk x_I>}e^{<\frac{\vk c^*_I}{u},\vk x_I>}e^{<\vk c^*_J,\vk x_J>}e^{-\frac{1}{2}(\vk x/\bar{\vk u})^\top\Sigma^{-1}(\vk x/\bar{\vk u})}\\
&\leq&e^{< \vk \lambda^{\varepsilon}_I,\vk x_I>}e^{<\vk c_J^*,\vk x_J>}e^{-\frac{1}{2}(\vk x/\bar{\vk u})^\top\left(\Sigma^{-1}-\varepsilon I\right)(\vk x/\bar{\vk u})}e^{-\frac{\varepsilon}{2}||\vk x/\bar{\vk u}||^2}\\
&\leq&e^{< \vk \lambda^{\varepsilon}_I,\vk x_I>}e^{<\vk c_J^*,\vk x_J>}e^{-\frac{\varepsilon}{2}||\vk x_J||^2}=:\bar\theta(\vk x),
}
where
\bqny{
\vk\lambda^\varepsilon:= \vk\lambda^\varepsilon(\vk x)=\vk \lambda +sign(\vk x)\varepsilon> \vk 0_I, \quad \vk c^*=\vk c^\top\Sigma^{-1}.
}

Let $\vk x\in\R^d$ and define two sets $Z_+,~Z_-\subset\{1\ldot d\}$ such that $\forall i\in Z_+~ x_i>0$, and $\forall i\in Z_-~ x_i\leq 0$. For all large $u$ we have
\bqny{
h_u([0,\CN{\Lambda}],\vk x)&=&\pk{\exists_{t\in[0,\CN{\Lambda}]}:u(\Gamma \vk B_u(t)+\overline{t}\vk v -\overline{t}\vk c-u\vk a)>\vk 0}\\
&\leq& \pk{\exists_{t\in[0,\CN{\Lambda}]}:u(\Gamma \vk B_u(t)+\overline{t}\vk v -\overline{t}\vk c-u\vk a)_{Z_+}>\vk 0_{Z_+}}\\
&\leq& \pk{\exists_{t\in[0,\CN{\Lambda}]}:u(\Gamma \vk B_u(t)+\overline{t}\vk v -\overline{t}\vk c-u\vk a)_{Z_+\cap I}>\vk 0_{Z_+\cap I}}\\
&=& \pk{\exists_{t\in[0,\CN{\Lambda}]}:u(\Gamma \vk B_u(t))_{Z_+\cap I}-(1-t/u^2)\vk x_{Z_+\cap I} - t\vk a_{Z_+\cap I}>\vk 0_{Z_+\cap I}}\\
&\leq&\pk{\exists_{t\in[0,\CN{\Lambda}]}:\sum_{i\in Z_+\cap I}u(\Gamma \vk B_u(t))_{i}>C_1\sum_{i\in Z_+\cap I}x_i + C_2},
}
where
\bqny{
u&>&2\sqrt{T}, \quad C_1=1-\frac{T}{4T}=\frac{3}{4}, \quad
C_2=-T\sum_{i\in I}|\vk a_i|.
}
 Moreover, using the properties of Brownian bridge we know that if $u>1$, then for any $s,t\in[0,T]$ and $i\in\{1\ldot d\}$ we have
\bqn{
u^2\E{(\Gamma \vk B_u(t)_i-\Gamma \vk B_u(s)_i)^2}&=&2u^2\left(\E{\left.(\vk W(\overline{t})_i-\vk W(\overline{s})_i)^2\right|\vk W(1)=\vk v}+\vk v^2_i(\overline{t}-\overline{s})^2\right)\\
&\leq& 2u^2\left(C_0|\overline{t}-\overline{s}|+\vk v^2_i(\overline{t}-\overline{s})^2\right)\notag\\
&\leq& 2\left(C_0|t-s|+\frac{\vk v^2_i}{u^2}(t-s)^2\right)\notag\\
&\leq& 2\left(C_0+T\vk v^2_i\right)|t-s|\notag
}
for some constant $C_0>0$. Using further the Piterbarg inequality we have for some constant $C>0$
\bqny{
h_u([0,\CN{\Lambda}],\vk x)&\leq& C \left(C_1\sum_{i\in Z_+\cap I}x_i + C_2\right)e^{-\frac{1}{2\sigma^2(\CN{\Lambda})}\left(C_1\sum_{i\in Z_+\cap I}x_i + C_2\right)^2}
\mathbb{I}{\left\{\sum\limits_{i\in Z_+\cap I}x_i>-C_b\right\}}+\mathbb{I}{\left\{\sum\limits_{i\in Z_+\cap I}x_i\leq -C_b\right\}}\\
&\leq& \bar{C}e^{-C^*\sum_{i\in Z_+\cap I}x_i^2}
\mathbb{I}{\left\{\sum\limits_{i\in Z_+\cap I}x_i>C_{b}\right\}}+\mathbb{I}{\left\{\sum\limits_{i\in Z_+\cap I}x_i\leq C_{b}\right\}}\\
&=:&\bar{h}([0,\CN{\Lambda}],\vk x),
}
where $C_b=-C_2/C_1$.
We may assume that $\bar{C}$ and $C^*$ are not dependent on $Z_+$.
For every non-empty set $Z\subset\{1\ldot d\}$
define $\mathbb{S}_{Z}=\{x\in\R^d:\forall i\in Z,~x_i>0,~\forall j\not\in Z,~x_j\leq 0\}$, where
\bqny{
A_Z=\left\{\vk x\in \mathbb{S}_Z: \sum\limits_{i\in Z\cap I}x_i\leq C_{b}\right\}.
}
Setting   $K=I \cap Z, L=I\setminus Z$ we may further write
\bqny{
\lefteqn{ \int_{\mathbb{S}_{Z}}\bar{h}([0,\CN{\Lambda}],\vk x)\bar{\theta}(\vk x)\td\vk x}\\
&\leq&\int_{A_Z}\bar{\theta}(\vk x)\td\vk x+\int_{\mathbb{S}_Z\setminus A_Z}\bar{h}([0,\CN{\Lambda}],\vk x)\bar{\theta}(\vk x)\td\vk x\\
&\leq&\int_{A_Z}\bar{\theta}(\vk x)\td\vk x+\bar{C}\int_{\mathbb{S}_Z\setminus A_Z}e^{-C^*\sum_{i\in Z_+\cap I}\vk x_i^2}e^{< \vk \lambda^{\varepsilon}_I,\vk x_I> +<\vk c_J^*\vk x_J> -\frac{\varepsilon}{2}||\vk x_J||^2}\td\vk x\\
&\leq&\int_{{(A_Z)}_J}e^{<\vk c_J^*,\vk x_J>-\frac{\varepsilon}{2}||\vk x_J||^2}\td\vk x_J\int_{{(A_Z)}_{ { K }}}e^{\sum_{i\in K}\vk\lambda^{\varepsilon}_i\vk x_i}\td\vk x_{ { K }} \int_{{(A_Z)}_{ { L }}}e^{\sum_{i\in  { L }}\vk\lambda^{\varepsilon}_i \vk x_i}\td\vk x_{ { L }}\\
&+&\bar{C}\int_{(\mathbb{S}_Z\setminus A_Z)_J}e^{<\vk c_J^*,\vk x_J>-\frac{\varepsilon}{2}||\vk x_J||^2}\td\vk x_J \int_{(\mathbb{S}_Z\setminus A_Z)_{ { K }}} e^{-\sum_{i\in  { K }}(C^*\vk x_i^2-\vk\lambda^{\varepsilon}_i\vk x_i) }\td\vk x_{ { K }}\\
& & \times\int_{(\mathbb{S}_Z\setminus A_Z)_{ { L }}}e^{\sum_{i\in  { L }}\vk\lambda^{\varepsilon}_i\vk x_i}\td\vk x_{ { L }}.
}
We show next that all the above  integrals are finite. We have
\bqny{
& &\int_{{(A_Z)}_J}e^{<\vk c_J^*,\vk x_J>-\frac{\varepsilon}{2}||\vk x_J||^2}\td\vk x_J\leq\int_{\R^{|J|}}e^{<\vk c_J^*,\vk x_J>-\frac{\varepsilon}{2}||\vk x_J||^2}\td\vk x_J<\infty,\\
& &
\int_{(\mathbb{S}_Z\setminus A_Z)_J}e^{<\vk c_J^*,\vk x_J>-\frac{\varepsilon}{2}||\vk x_J||^2}
d\vk x_J\leq \int_{\R^{|J|}}e^{<\vk c_J^*,\vk x_J>-\frac{\varepsilon}{2}||\vk x_J||^2} d\vk x_J<\infty,\\
& &\int_{(\mathbb{S}_Z\setminus A_Z)_{ { K }}} e^{-\sum_{i\in  { K }}(C^*\vk x_i^2-\vk\lambda^{\varepsilon}_i\vk x_i)}\td\vk x_{ { K }}\leq \int_{\R^{| { K }|}} e^{-\sum_{i\in  { K }}(C^*\vk x_i^2-\vk\lambda^{\varepsilon}_i\vk x_i)}\td\vk x_{ { K }}<\infty
}
and for the other two integrals (recall that $\vk\lambda^{\varepsilon}_I$ has positive components)
\bqny{
& &\int_{{(A_Z)}_{I\setminus Z}}e^{\sum_{i\in I\setminus Z}\vk\lambda^{\varepsilon}_i\vk x_i}\td\vk x_{ { L }}\leq\int_{(-\IF, 0)^{| { L }|}}e^{\sum_{i\in  { L }}\vk\lambda^{\varepsilon}_i\vk x_i}\td\vk x_{ { L }}<\infty,\\
& &\int_{(\mathbb{S}_Z\setminus A_Z)_{ { L }}}e^{\sum_{i\in  { L }}\vk\lambda^{\varepsilon}_i\vk x_i}\td\vk x_{ { L }}\leq \int_{(-\IF, 0)^{| { L }|}}e^{\sum_{i\in  { L }}\vk\lambda^{\varepsilon}_i\vk x_i}\td\vk x_{L}<\infty,
}
where the last term is finite because the area $(A_{Z})_{ { K }}$ is finite.

Consequently, the dominated convergence theorem yields
\bqn{
M(u,\CN{\Lambda})\sim F(u)\int_{\R^d}h([0,\CN{\Lambda}],\vk x)\theta(\vk x)\td\vk x
,\qquad u\to\infty.\label{Mint}
}
 Suppose for simplicity that the number of elements of $J$ is $d-m>0$ and define
\bqny{
\widehat{C}&=&e^{\frac{1}{2}\tilde{\vk c}^\top(\Sigma^{-1})_{JJ}\tilde{\vk c}}, \ \
\nu_{J}(\vk x_J)=e^{-\frac{1}{2}(\vk x_J-\tilde{\vk c})^\top(\Sigma^{-1})_{JJ}(\vk x_J-\tilde{\vk c})},\\
Q&=&\int_{\R^{d-m}}\mathbb{I}_{\{\vk x_{U}<\vk 0_{U}\}}\nu_J(\vk x_J)\td\vk x_J\in(0,\infty),\\
E([0,\CN{\Lambda}])&=&\int_{\R^m}\pk{\exists_{t\in[0,\CN{\Lambda}]}:\vk W_{I}(t)-t\vk a_{I}>\vk x_I}e^{\vk \lambda_I^\top \vk x_I}\td\vk x_I.
}
Then we can rewrite \eqref{Mint} as
\bqny{
M(u,\CN{\Lambda})\sim \widehat{C} Q E([0,\CN{\Lambda}])F(u),\qquad u\to\infty.
}
Since $\tilde{\vk c}=\vk c_J+((\Sigma^{-1})_{JJ})^{-1}(\Sigma_{JI})^{-1}\vk c_I$, then
$$ \vk c^\top\Sigma^{-1}\vk c= \vk c_I^\top(\Sigma_{II})^{-1}\vk c_I+ 2\log\widehat{C}$$
implying that
\bqny{
\widehat{C}F(u)=\frac{(2\pi)^{|I|/2}\sqrt{|\Sigma_{II}|}}{(2\pi)^{d/2}\sqrt{|\Sigma|}}u^{-|I|}\varphi_I(u\vk a_I+\vk c_I),
}
where $\varphi_I$ is pdf of  $\vk W_{I}(1)$. Since for  $\vk\xi_{J}$  a Gaussian random vector with mean
 $\tilde{\vk c}$ and covariance matrix $((\Sigma^{-1})_{JJ})^{-1}$ its  pdf is
\bqny{
\frac{1}{(2\pi)^{\kk{|J|/2}}\sqrt{|\Sigma_{JJ}|}}\nu_{J}(\vk x_J), \ \ \vk x\inr^d
}
and thus
\bqny{
Q=(2\pi)^{|J|/2}\sqrt{|\Sigma_{JJ}|}\pk{\vk \xi_U>\vk c_U},
}
which implies
\bqny{
\widehat{C} Q F(u)=\pk{ \vk \xi_U> \vk c_U}u^{-|I|}\varphi_I(u\vk a_I+\vk c_I),
}
where $\vk \xi_U$ is $|U|$-dimensional Gaussian vector with expectation $-(((\Sigma^{-1})_{JJ})^{-1})_{AJ}(\Sigma^{-1})_{JI}\vk c_I$ and covariance matrix $(((\Sigma^{-1})_{JJ})^{-1})_{UU}$. Hence the proof follows since $ \vk \xi_U$ has the same distribution as
${\vk W_U(1) \lvert \vk W_I= \vk c_I}$. We show next the claim in \eqref{T}.
For any $\CN{\Lambda}>0$, since by \eqref{MC} $\lambda_i>0$ for all $i\in I$, we have
\bqny{
E([0,\CN{\Lambda}])&\leq&\int_{\R^m}\pk{\exists_{t\in[0,\infty]}:\vk W_{I}(t)-t\vk a_{I}>\vk x_I}e^{\vk \lambda_I^\top \vk x_I}\td\vk x_I\\
&\leq&\sum_{i=0}^{\infty}\int_{\R^m}\pk{\exists_{t\in[i,i+1]}:\vk W_{I}(t)-t\vk a_{I}>\vk x_I}e^{\vk \lambda_I^\top \vk x_I}\td\vk x_I\\
&\leq&\sum_{i=0}^{\infty}\int_{\R^m}\pk{\sup\limits_{t\in[i,i+1]}(\vk W_{I}(t)-t\vk a_{I})>\vk x_I}e^{\vk \lambda_I^\top \vk x_I}\td\vk x_I\\
&=&\frac{1}{\prod\limits_{i\in I}\lambda_i}\sum_{i=0}^{\infty}\E{e^{\sprod{\vk\lambda_I,\vk M_i}}},
}
where $m$ is the number of elements of $I$ and
\bqny{
\vk M_i=\sup_{t\in[i,i+1]}(\vk W_{I}(t)-t\vk a_{I}).
}
Using the independence of increments of Brownian motion we have
\bqny{
\vk M_i=\sup_{t\in[0,1]}(\vk W_{I}(t)-t\vk a_{I})+\vk V_I(i) - i\vk a_I,
}
where $\vk V_I$ is an independent copy of $\vk W_I$. By the aforementioned independence, we have
\bqny{
\E{e^{\sprod{\vk\lambda_I,\vk M_i}}}=\E{e^{\sprod{\vk\lambda_I,\sup_{t\in[0,1]}(\vk W_{I}(t)-t\vk a_{I})}}}\E{e^{\sprod{\vk\lambda_I,\vk V_I(i)}}}e^{-i\sprod{\vk\lambda_I,\vk a_I}}.
}

Since
\bqny{
\sprod{\vk\lambda_I,\vk V_I(i)}=\sprod{\vk \lambda,\Gamma \vk B(i)}=\sprod{\Gamma^\top\vk\lambda,\vk B(i)}
}
we obtain further
\bqny{
\E{e^{\sprod{\vk\lambda_I,\vk V_I(i)}}}
=\prod_{j=1}^{d}e^{\frac{i(\Gamma ^\top\vk\lambda)_j^2}{2}}
=e^{\frac{i \vk\lambda^\top\Sigma\vk\lambda}{2}}=e^{\frac{i \tilde{\vk a}^\top \Sigma^{-1} \tilde{\vk a}}{2}}
}
and consequently
\bqny{
\frac{1}{\prod\limits_{i\in I}\lambda_i}\sum_{i=0}^{\infty}\E{e^{\sprod{\vk\lambda_I,\vk M_i}}}&=&\E{e^{\sprod{\vk\lambda_I,\sup_{t\in[0,1]}(\vk W_{I}(t)-t\vk a_{I})}}}\sum_{i=0}^{\infty}e^{-i\sprod{\vk\lambda_I,\vk a_I}}e^{\frac{i \tilde{\vk a}^\top \Sigma^{-1} \tilde{\vk a}}{2}}\\
&=&\E{e^{\sprod{\vk\lambda_I,\sup_{t\in[0,1]}(\vk W_{I}(t)-t\vk a_{I})}}}\sum_{i=0}^{\infty}e^{\frac{-i \tilde{\vk a}^\top \Sigma^{-1} \tilde{\vk a}}{2}}
}
since by \eqref{MC} and \eqref{eq:IJ3} $\sprod{\vk\lambda_I,\vk a_I}=\vk a^\top_I (\Sigma_{II})^{-1}\vk a_I= \tilde{\vk a}^\top \Sigma^{-1} \tilde{\vk a}>0$.
Hence the proof follows from the monotone convergence theorem.
\QED

\section*{Acknowledgments}
We are thankful to the reviewers for valuable comments and corrections.
K.D. was partially supported by NCN Grant No
2018/31/B/ST1/00370.
Thanks to the Swiss National Science Foundation Grant 200021-196888.

\bibliographystyle{ieeetr}
\bibliography{EEEA}

\end{document}